\input amstex.tex
\documentstyle{amsppt}
\magnification=1200
\hsize=150truemm
\vsize=224.4truemm
\hoffset=4.8truemm
\voffset=12truemm

\TagsOnRight
\NoBlackBoxes
\NoRunningHeads

\def\Square{\rlap{$\sqcup$}$\sqcap$}
\def\cqfd {\quad \hglue 7pt\par\vskip-\baselineskip\vskip-\parskip
{\rightline{\Square}}}

\define\C{{\Cal C}}
\define\F{{\Cal F}}
\define\A{{\Cal A}}

\redefine\P{{\Cal P}}
\define\G{{\Cal G}}
\define\R{{\bold R}}
\define\Z{{\bold Z}}
\define\N{{\bold N}}
\define\pim{\pi _1M}
\define\Rt {${\bold R}$-tree}
\define\mi{^{-1}}
\define\st   {\text{\rm Stab}\,}
\define\rk{\text{\rm rk}\,}
\define\fix{\text{\rm Fix}\,}
\let\thm\proclaim
\let\fthm\endproclaim
\let\inc\subset 
\let\ds\displaystyle
\let\wtilde\widetilde 
\let\ov\overline

\define\fp{fixed point}
\let\bo\partial
\define\Aut{\text{\rm Aut}\,}
\define\Out{\text{\rm Out}\,}
\define\da{\partial  \alpha }
\define\db{\partial  \beta }
\define\dg{\partial  \Gamma  }
\define\hr{H\"older}
\define\bscp#1#2{\bigl(#1|#2\bigr)}
\define\scp#1#2{(#1|#2)}
\define\tg{{\wtilde{\Cal G}}}
\define\h{{\Cal H}}
\define\th{{\wtilde{\Cal H}}}

\define\PF{_{PF}}
\define\dt{I}
\define\pt{I}
\define\pth{II}

\newcount\tagno
\newcount\secno
\newcount\subsecno
\newcount\stno
\global\subsecno=1
\global\tagno=0
\define\ntag{\global\advance\tagno by 1\tag{\the\tagno}}

\define\sta{\the\secno.\the\stno
\global\advance\stno by 1}

\define\stap{\the\secno.\the\stno
\global\advance\stno by 1} 

\define\sect{\global\advance\secno by 1\global\subsecno=0\global\stno=1\
\the\secno. }

\define\subsect{}

\def\nom#1{\edef#1{\the\secno.\the\stno}}
\def\snom#1{\edef#1{\the\secno.\the\stno}}
\def\eqnom#1{\edef#1{(\the\tagno)}}

\newcount\refno
\global\refno=0
\def\nextref#1{\global\advance\refno by 1\xdef#1{\the\refno}}
\def\bref {\ref\global\advance\refno by 1\key{\the\refno}}

\def\strutdepth{\dp\strutbox}
\def \ss{\strut\vadjust{\kern-\strutdepth \sss}} 
\def \sss{\vtop to \strutdepth{
\baselineskip\strutdepth\vss\llap{$\diamondsuit\;\;$}\null}}

\def\strutdepth{\dp\strutbox}
\def \sst{\strut\vadjust{\kern-\strutdepth \ssss}}
\def \ssss{\vtop to \strutdepth{
\baselineskip\strutdepth\vss\llap{$\spadesuit\;\;$}\null}}

\nextref\Be
\nextref\BF
\nextref\BFHg
\nextref\BFHt
\nextref\BFHk
\nextref\BFHs
\nextref\BH
\nextref\Bow
\nextref\BHa
\nextref\CL
\nextref\CLvs
\nextref\CT
\nextref\Dah
\nextref\DV
\nextref\Far
\nextref\Gab
\nextref\GJLL
\nextref\GL
\nextref\Guia
\nextref\Gui
\nextref\HW
\nextref\Hi
\nextref\Ka
\nextref\KW
\nextref\Lev
\nextref\LLun
\nextref\LLdeux
\nextref\LLquat
\nextref\LLdeuxdeux
\nextref\Ludis
\nextref\Luconj
\nextref\Mi
\nextref\Ne
\nextref\Pau
\nextref\Pa
\nextref\Se
\nextref\Sh
\nextref\Swa

\topmatter

\abstract
We show that every automorphism $\alpha$ of a free group $F_k$ of finite
 rank $k$ 
has {\it asymptotically periodic\/} dynamics on $F_k$ and its
boundary $\partial F_k$: there exists a positive power $\alpha^q$ such
that every element of the compactum $F_k \cup \partial F_k$ converges to a
fixed point under iteration of $\alpha^q$.

Further results about the dynamics of $\alpha$, as well as an  
extension from $F_k$ to word-hyperbolic groups, are given in the later
sections.

\endabstract

\title Automorphisms of free groups have asymptotically 
periodic dynamics
    \endtitle

\author  Gilbert Levitt, Martin Lustig 
\endauthor

\toc
\widestnumber\head{10}
\specialhead {} \, Introduction and statement of results\endspecialhead
\head 1. Preliminaries\endhead
\subhead 1.a. $F_k$ and its boundary\endsubhead
\subhead 1.b. Limit sets and asymptotic periodicity\endsubhead
\subhead 1.c. Invariant trees\endsubhead
\subhead 1.d. Bounded backtracking and $Q(X)$\endsubhead
\head 2. A lemma on asymptotic periodicity\endhead
\head 3. Limit sets of interior points\endhead
\head 4. The simplicial case\endhead
\head 5. A reduction\endhead
\head 6. The train track and the tree\endhead
\subhead {} \ \ Creating a shortcut\endsubhead
\subhead {} \ \  Legal paths\endsubhead
\subhead {} \ \  The PF-metric and the invariant tree\endsubhead
\subhead {} \ \  Elliptic elements\endsubhead
\head 7. Geometry on the train track\endhead
\subhead {} \ \  Spaces quasi-isometric to trees\endsubhead
\subhead {} \ \  $K$-PF-geodesics\endsubhead
\subhead {} \ \  An inequality\endsubhead
\head 8. Proof of Theorem 5.1\endhead
\subhead {} \ \  Paired train tracks\endsubhead
\subhead {} \ \  Creating a fixed point\endsubhead
\subhead {} \ \  The main argument\endsubhead
\newpage
\head 9. More on the dynamics\endhead
\subhead {} \ \  Products of trees\endsubhead
\subhead {} \ \  Dynamics of irreducible automorphisms\endsubhead
\subhead {} \ \  The number of periods\endsubhead
\subhead {} \ \  Automorphisms with many fixed points\endsubhead
\head 10. Hyperbolic groups\endhead
\subhead {} \ \  Bounding periods\endsubhead
\subhead {} \ \  One-ended groups\endsubhead
\subhead {} \ \  Free products\endsubhead
\subhead {} \ \  Groups with torsion\endsubhead
\head 11. Examples and questions\endhead
\subhead {} \ \  Examples\endsubhead
\subhead {} \ \  Free groups\endsubhead
\subhead {} \ \  Hyperbolic groups\endsubhead
\subhead {} \ \  Actions with finite limit sets\endsubhead
\head {} References\endhead
\endtoc

\endtopmatter

\document

\head  {Introduction and statement of results}\endhead

Let
$F_{k}$ denote the free group of   rank $k \geq 2$.
Conjugation $i_u$ by an element $u\in F_k$    has very simple
dynamics. If $g\in F_k$
     commutes with $u$, then $g$  is a fixed point of $i_u$. If $g$ does
not commute with $u$, then the length of  $(i_u)^n(g)=u^ngu^{-n}$
tends to infinity as $n\to+\infty$, and $u^ngu^{-n}$ converges to the
infinite word
$u^\infty=\lim_{n\to+\infty}u^n$.

On the space of infinite words (which may be viewed as the boundary
$\bo F_k$),  the action of $i_u$ is simply
left-translation by $u$. It has North-South (loxodromic) dynamics:
$u^\infty$ is an attracting fixed point (a sink), $u^{-\infty}$ is a
repelling fixed
point (a source), and $\lim_{n\to\pm\infty}u^nX=u^{\pm\infty}$ for
every infinite word
$X\neq u^{\pm\infty}$.
Similar considerations apply to conjugation by any element of infinite
order in a word hyperbolic group.

We proved in [\LLun] that ``most automorphisms'' (in a precise sense)
of a given
hyperbolic group (e.g\. $F_k$)
    have   North-South dynamics on the
boundary of the group. But, of course, interesting automorphisms
usually are not ``generic''.

For
instance, Nielsen studied mapping classes of surfaces by lifting them to the
universal covering, and considering the action of various lifts
    on the circle at infinity  $S_\infty$ (the boundary of the surface group).
He used lifts with more than two periodic points (see [\Mi]),
and  one of his key
results is that   a lift
$f$ always has periodic points on  $S_\infty$ (equivalently, its
rotation number on  $S_\infty$ is rational).

Since $f$ induces a homeomorphism $\bo f$ of the circle  $S_\infty$,
there is a power $h$ of $\bo f$ such that, for  
any
$X$ on the circle,   the sequence
$h^n(X)$ converges to a fixed point of $h$ as $n\to+\infty$;  we will say
that $\bo f$ has   asymptotically periodic dynamics.


Our main results may be viewed as a generalization of these facts to
arbitrary automorphisms of free (or hyperbolic) groups.

Let
$\alpha $ be an automorphism of $F_k$.
It  induces canonically a homeomorphism $\partial
\alpha$ on the   boundary $\partial F_{k}$.  The latter, a Cantor
set, can be identified with the set of reduced right-infinite words in
an arbitrary basis of $F_{k}$, or with the space of ends of any
simplicial tree on which $F_{k}$ acts freely.  As usual, we provide
$F_{k}$ with the discrete topology and use $\partial F_{k}$ to
compactify $F_{k}$, thus obtaining $\ov F_{k} = F_{k} \cup \partial
F_{k}$.  One obtains from $\alpha$ and $\partial \alpha$ together a
homeomorphism $\ov \alpha: \ov F_{k} \to \ov F_{k}$ of this
compactum.  This paper studies the dynamics of this homeomorphism.


The automorphism
$\alpha
\in \Aut (F_{k})$ has {\it asymptotically periodic dynamics\/} (on $\ov F_k$)
if there exists $q\ge1$ such that, for every $X
\in
\ov F_{k}$, the sequence
$\ov\alpha ^{qn}(X)$ converges (to a fixed point of $\ov\alpha ^{q}$).
Our main
result can now be stated
as follows:

\thm
{  Theorem \pt} Every automorphism $\alpha \in \Aut (F_{k})$ has
asymptotically periodic
dynamics on  $\ov F_k$.
\fthm

It follows from [\LLdeux] that $q$ may be bounded by a number $M_k$ depending
only on $k$ (see Theorem 1.1 below).

If $g\in F_k$ is not $\alpha $-periodic, then by Theorem I   the set
of limit
points of the sequence
$\alpha ^n(g)$ as $n\to+\infty$ is a periodic orbit of $\da$. In
other words, there
exists a positive  integer   $q$ such that, for any $N$, the sequence
consisting of the initial segment of length $N$ of $\alpha ^n(g)$ is
eventually periodic with period dividing $q$.

As an illustration, define an   automorphism   on
the free group of rank 3 by $\alpha (a)=cb$, $\alpha (b)=a$,  $\alpha
(c)=ba$. Applying powers of $\alpha $ to $a$ gives $a\mapsto
cb\mapsto baa\mapsto acbcb\mapsto cbbaabaa\mapsto
baaacbcbacbcb\mapsto\dots$, showing that $\alpha ^n(a)$ limits onto
an orbit of period 3. On the other hand $a\mi\mapsto b\mi
c\mi\mapsto a\mi a\mi b\mi\mapsto b\mi c\mi b\mi c\mi
a\mi\mapsto\dots$, and
$\alpha ^n(a\mi)$ limits onto an orbit of period 2. 

Other examples
will be given in
section 11.  In particular, it  is quite common for a boundary point 
of an $\alpha
$-invariant free factor $F$ of $F_k$   to be the limit of orbits in  $\ov
F_k$ well outside
$\ov F$.

For an orientation-preserving homeomorphism of the circle, all
periodic points have the same period, and existence of a periodic
point is enough to imply asymptotically periodic dynamics. For
automorphisms of free groups, there may be periodic points with
different periods on the boundary (as shown by the above example). It
is relatively
easy to prove that periodic points exist, but much harder to control all
sequences $\ov\alpha ^{n}(X)$.


\medskip

More generally, there are many examples of actions of a group $\Gamma $ on a
compact space $X$ with the
following property: {\it There exists $q\ge1$ such that, for every $x\in X$ and
$g\in \Gamma $, the sequence $g^{qn}(x)$ converges as $n\to+\infty$.}
For instance:

$\bullet$ The action of $\Aut(F_k)$ on $\bo F_k$ and $\ov F_k$.
Conjecturally, the
action of
$\Aut (\Gamma )$ on $\ov \Gamma $, for an arbitrary hyperbolic group 
$\Gamma $ (this is true if $  \Gamma $ is one-ended and virtually
torsion-free, see Theorem \pth{} below).

$\bullet$ Convergence actions of virtually torsion-free groups.

$\bullet$ The action of the mapping class group of a closed surface on the
Thurston boundary of Teichm\"uller space (this follows from Nielsen-Thurston
theory). By analogy, one may ask about the action of $\Out (F_k)$ on
the boundary
of Culler-Vogtmann's outer space (see [\BFHk], [\CLvs],  [\LLquat]
for   partial
results).

$\bullet$ The action of $\pim$ on the sphere at infinity $S_\infty$ of
$\wtilde M$, where
    $  M$ is a closed Riemannian manifold  (or orbifold) with negative
curvature and
$\wtilde M$ is the universal covering. Flat manifolds also provide
  examples, because
of Bieberbach's theorem. This does not extend  to
arbitrary non-positively
curved manifolds, but the above property might hold on an open dense subset
of $X=S_\infty$.
 

Knowing that a group $\Gamma $ acts on $X$ with the above property
gives a lot of
information about dynamics of individual elements of $\Gamma $. It also
gives the   global information   that $q$ does not
depend on
$g$. There may exist a  stronger property, that would be weaker than the
convergence property but strong enough to contain more global information on
$\Gamma $. Let us also mention the Powers   property, which implies
$C^*$-simplicity (see [\BHa]). We will observe in \S\kern.15em 11 that, if
$\Gamma $ is a  non-elementary torsion-free hyperbolic group, and 
$H$ is a      subgroup  of $\Aut(\Gamma )$ containing all inner
automorphisms, then the action of $H$   on
$\bo\Gamma
$ has the Powers property. 

\medskip

Returning to automorphisms of free groups, the basic tool in our approach 
is  an {\it
$\alpha$-invariant}
$\R$-tree, i.e\. an $\R$-tree $T$ with an action of $F_{k}$ by isometries
which is minimal, non-trivial, with {\it trivial arc stabilizers\/}, and
$\alpha$-invariant: its length function $\ell$ satisfies
$\ell\circ\alpha =\lambda
\ell $ for some $\lambda \ge1$. The automorphism $\alpha $ is then
realized on $T$, in
the sense that there is
    a homothety
$H: T
\to T$, with stretching factor $\lambda \geq 1$, such that
    $\alpha(w) H = H w: T \to T$ for all $w \in F_{k}$. If
$\lambda = 1$, the tree  $T$ is simplicial and $H$ is an isometry. If
$\lambda > 1$,
    the action of $F_{k}$ on $T$ is non-discrete; in fact, every
$F_{k}$-orbit is dense in $T$.

In most cases, the map $H$ has a fixed point $Q$ (in $T$ or in its metric
completion $\ov T$).  The stabilizer of   $Q  $   is an
$\alpha$-invariant subgroup $\st Q\inc F_{k}$, which has rank
    strictly smaller than $k$
    [\GL]. This allows us to set up the proof of our main result as a
proof by induction over the rank $k$.

For $g\in F_k$, we study the behavior of the sequence $\alpha ^n(g)$
through that of the sequence $\alpha ^n(g)Q=H^n(gQ) $, where $Q$ is a
fixed point of
$H$.

There are three main cases (see \S\kern.15em 3).

- If $gQ=Q$, we use the induction hypothesis.

- If $\lambda >1$ and $H^n(gQ) $ goes out to infinity in $T$ in a
definite direction,
then  $\alpha ^n(g)$ converges to an attracting fixed point of $\da$.

- If $H^n(gQ) $ ``turns around $Q$'', one shows that $\alpha ^n(g)$
accumulates onto
a periodic orbit contained in $\bo\st Q$, using a cancellation argument given
in \S\kern.15em 2.

Similar arguments (given in \S\kern.15em 4) make it possible to understand the
behavior of $\da^n(X)$, for $X\in\bo F_k$, when there   are enough  
simplicial invariant
trees, in particular when $\alpha $ is a  polynomially growing automorphism.

The general case is dealt with in \S\S\kern.15em 5 through 8. It
makes use of the point
$Q(X)$ introduced in [\LLquat], 
which reflects the dynamics of $\da$ on $X$ and allows us to extend
the approach from the  special cases dealt with previously.
The proof consists of geometric arguments on
relative train
tracks, and involves the 
asymptotic behavior of four distinct ways of measuring length  under iteration
of $\alpha$ (and of $\alpha^{-1}$). A sketch of the proof will be given at
the beginning of \S\kern.15em 5.

In \S\kern.15em 9, we prove a few more results about dynamics of
automorphisms of
free groups. We show that $F_k$ acts discretely on a suitable product
of trees (a
result proved in [\BFHg] and [\Ludis] for irreducible automorphisms).
We study the
bipartite graph whose vertices are the  attracting and repelling
fixed points of
$\da$, for
$\alpha $ irreducible. We show that, for an arbitrary $\alpha $, the
number of different
periods appearing in the dynamics of
$\ov\alpha $ is bounded by a number depending only on $k$ and
growing roughly like
$e^{\sqrt k}$. We also give a short proof of a result of [\BFHs] constructing
automorphisms with many fixed points.

In \S\kern.15em 10, we explain how to adapt the  arguments of
\S\S\kern.15em 2, 3, 4
     to a hyperbolic group $\Gamma $.   With asymptotically
periodic dynamics defined as in the case of free groups, we show: 

\thm
{  Theorem \pth }  Let $\alpha \in\Aut(\Gamma )$, with $\Gamma $ a virtually
torsion-free
hyperbolic group. \roster
\item Periodic  orbits of   $\ov \alpha $ have  at most
$M $ points, with $M$ depending only on $\Gamma $.
\item Every $g\in \Gamma $ is asymptotically periodic: there exists a
positive integer $q$ (depending   only on $\Gamma $) such that
$\ov\alpha^{qn}(g)
$ converges.
\item If $\Gamma  $ is one-ended, or $\alpha $ is polynomially
growing, then $\alpha $ has asymptotically periodic dynamics on $\ov\Gamma $.
\endroster
\fthm

It is not known whether all hyperbolic group are virtually
torsion-free (see [\KW]). In any case, it
    seems reasonable to conjecture that all automorphisms of
hyperbolic groups have asymptotically periodic dynamics. 

We start  \S\kern.15em 10 by  giving a proof of an unpublished result by
Shor [\Sh]: Up to
isomorphism, there are only finitely many fixed subgroups in a given
torsion-free
hyperbolic group. As in [\Sh], we use results by Sela [\Se], Guirardel [\Guia],
Collins-Turner [\CT], but we simplify the proof by using the fact
(proved in [\Lev]) that
$\Aut(\Gamma )$ contains only finitely many torsion conjugacy
classes, for $\Gamma $ a
torsion-free hyperbolic group.

Theorem \pth{} is first proved for torsion-free groups. For
one-ended groups, we use the
simplicial tree (with
cyclic edge stabilizers) given by the JSJ splitting. For
free products, we use  an \Rt{} constructed from the train tracks of [\CT].
Finally, we extend
our results to virtually torsion-free groups.

We conclude the paper by a section devoted to  examples and questions.

\bigskip

\head  {\sect Preliminaries}\endhead

\subhead \subsect  \the\secno.a. $F_k$ and its boundary\endsubhead

Let $F_k$ be a free group of rank $k\ge2$. Its boundary (or
space of ends) $\bo F_k$    is a Cantor set, upon which
$F_k$ acts by left translations. It compactifies $F_k$ into $\ov
F_k=F_k\cup\bo F_k$. The boundary $\bo J$ of a
    finitely generated subgroup
$J\inc F_k$ embeds naturally into $\bo F_k$. If $g\in F_k$ is
nontrivial, we let
$g^{\pm\infty}\in\bo F_k$ be the limit of $g^n$ as $n\to\pm\infty$.

If we choose a free basis $\A$ of $F_k$, we may view $\bo F_k$ as the  set of
right-infinite reduced words. The Gromov  scalar product $\scp  X
Y$ of two elements
$X,Y\in\ov F_k$ is the length of their maximal common initial
subword. A sequence
$X_n$ in $\ov F_k$  converges to $X\in\bo F_k$ if and only if $\scp
{X_n}X\to\infty$.

An
automorphism $\alpha \in\Aut (F_k)$ is a quasi-isometry of $F_k$. It
    induces a homeomorphism $\da:\bo F_k\to\bo F_k$, and also a homeomorphism
$\ov\alpha =\alpha
\cup\da$ of the compact space $\ov F_k$.  The conjugation $g\mapsto
wgw\mi$ will be denoted by $i_w$. Note that $\bo i_w$ is
left-translation by $w$.

Let $E$ be a set, and $f:E\to E$ a bijection. 
An element 
$x\in E$ is {\it periodic\/}, with period $q\ge1$, if $f^q(x)=x$ and
$q$ is the smallest positive integer with this property. The set
$\{x,f(x),\dots, f^{q-1}(x)\}$ is a {\it periodic orbit\/} of order $q$. 

We recall   a few  known results. The existence of a uniform bound
$M_k$ is the only fact used later. 

\nom\rappel
\thm{Theorem \sta{} 
[\GJLL,   \LLdeux,
\LLdeuxdeux]} Let $\alpha \in \Aut (F_{k})$.
\roster
\item Every periodic orbit of $\ov\alpha $ has order bounded by
$M_k$, where $M_k$
depends only on $k$, and $M_{k} \sim {({k \, {\hbox {\rm
log}}(k)})^{1/2}}$ as $k\to\infty$.

\item $\da$ has at least two periodic points of period $\le 2k$.

\item A fixed point of $\da$ which
does not belong to the boundary of  the fixed subgroup $\fix\alpha $
is either attracting
or repelling (sink or source). The number $a(\alpha )$ of orbits of
the action of
$\fix\alpha
$ on the set of attracting fixed points satisfies  the
``index inequality" $\ds\ \rk\fix\alpha +\frac12a(\alpha )\le k$.

\item Every attracting fixed point of $\da$   is
superattracting with respect to the canonical H\"older structure of
$\partial F_{k}$, with attraction rate $\lambda_{i} \ge 1$. If
$\lambda _i>1$, then
$\lambda _i$ is the exponential growth rate of some conjugacy class
under iteration
of $\alpha $. There are at most
${3k- 2 \over 4}$ distinct such growth rates.

\endroster
\fthm

\subhead \subsect \the\secno.b. Limit sets   and asymptotic
periodicity\endsubhead

Let $f$ be a homeomorphism of a compact space $K$ (for instance $\bo
F_k$ or $\ov
F_k$). 
Given
$y\in K$, the {\it $\omega $@-limit set\/} $\omega (y,f)$, or simply
$\omega (y)$,
    is the set of limit points of the sequence $f^n(y)$ as $n\to +\infty$. It is
compact, and invariant under
$f$ and
$f\mi$. We observe:

\snom\fini
\thm {Lemma \stap}  Let $f$ be a homeomorphism of a compact space
$K$. Given $y\in
K$ and $q\ge1$, the following conditions  are equivalent:
\roster\item  $\omega (y)$ is finite, and has $q$ elements.
\item $\omega (y)$ is a periodic orbit of order $q$.
\item The sequence $f^{qn}(y)$ converges as $n\to +\infty$, and $q$
is minimal for
this property.

\endroster
Given $p\ge2$, the set $\omega (y,f^p)$ is finite if
and only if $\omega (y,f)$ is finite. \cqfd
\fthm

If these equivalent conditions hold, we say that  $y$ is {\it asymptotically
periodic.}   In particular, we have defined asymptotically periodic elements of
$\ov F_k$ (with respect to $\ov\alpha $).

We say that   $\alpha \in\Aut(F_k)$ has {\it asymptotically
periodic dynamics\/}   if every $X\in\ov F_k$ is
asymptotically
periodic. By assertion (1) of Theorem \dt, the period $q$ is bounded
independently of $X$, and this is equivalent to the definition given in the
introduction.


\subhead \subsect \the\secno.c. Invariant trees\endsubhead

As in [\LLdeux], we will use as our basic tool  an  $\alpha $@-invariant
\Rt{}
$T$ with trivial arc stabilizers. We summarize its main
properties.

\snom\arbre
\thm{Theorem \stap{} [\GJLL, \GL]} Given $\alpha \in\Aut(F_k)$, there exists 
an
\Rt{}
$T$ such that:
\roster
\item"(a)" $F_k$ acts on $T$ isometrically,  non-trivially, minimally, with
trivial arc  stabilizers.
\item"(b)" There exist   a number 
$\lambda
\ge1$ and a homothety
$H\colon\ T\to T$ with stretching factor $\lambda $ such that
$$
\alpha (g)H=Hg
$$
for all $g\in F_k$ (viewing elements of $F_k$ as isometries of $T$). 
\item"(c)"   
    If  $T$ is not simplicial, then $\lambda >1$ (and all $F_k$-orbits are 
dense).
\item"(d)"
    Given $Q\in T$, its stabilizer $\st Q$ has rank $\le
k-1$, and the action of $\st Q$ on $\pi _0(T\setminus\{Q\})$ has at most $2k$
orbits. The number of $F_k$@-orbits of branch points of $T$ is at most $2k-2$.
\cqfd

\endroster\fthm

A tree with these properties will simply be called an {\it $\alpha $-invariant
\Rt\/}  (with trivial arc stabilizers). The pair $(\lambda ,H)$ is unique. A
construction will be sketched in \S \kern.15em 6.   The  length function of
$T$ is denoted by
$\ell:F_k\to[0,\infty)$, it satisfies
$\ell\circ\alpha =\lambda \ell$. An element
$g\in F_k$ is {\it elliptic\/} if $\ell(g)=0$ (i.e\. if $g$ has a
fixed point), {\it
hyperbolic\/} if
$\ell(g)>0$.

If a point $Q\in T$ with  nontrivial stabilizer  is fixed by $H$,
    the subgroup
$\st Q$ is
$\alpha $@-invariant. Since it has rank less that $k$, and $\partial  \st Q$
embeds into $\partial  F_k$, this will allow us  to use induction on $k$.

A {\it  ray\/} is (the image of) an isometric map $\rho $ from
$[0,\infty)$ or $(0,\infty)$ to $T$. It is an
{\it eigenray\/} of $H$ if
$\rho (\lambda t)=H\rho (t)$, a  {\it  periodic ray\/} if it is an
eigenray of some
power of $H$. As usual, the {\it boundary\/}
$\bo T$ is the set of equivalence classes of rays. The action of
$F_k$, and the map $H$,
extend to $\bo T$.

We also consider the action of $F_k$ on the metric
completion
$\ov T$ of
$T$. When $\lambda =1$, the   tree $T $ is simplicial, so $\ov T=T$.  
Now suppose $\lambda >1$.  Points
of
$\ov T\setminus T$ have trivial stabilizer.   
The homothety $H$ has
a  canonical extension to
$\ov T$, with a unique \fp{}
$Q\in
\ov T$. All  eigenrays have origin $Q$.
    If a component
     of $\ov T\setminus \{Q\}$  is fixed by $H$ (in particular if
$Q\notin T$), that component contains a unique   eigenray.

\subhead \subsect \the\secno.d. Bounded backtracking and $Q(X)$\endsubhead

Let $T$ be an $\alpha $-invariant \Rt{}  with trivial  arc stabilizers. Fix
$Q\in
\ov T$. When $\lambda >1$, we always choose
$Q$ to be the fixed point of $H$.

Let $Z$ be a geodesic metric space. We say that a map $f:Z\to
\ov T$ has {\it bounded backtracking\/} if there exists $C>0$ such
that the image of
any geodesic segment   $[P,P' ]$ is contained in the
$C$@-neighborhood of the segment $[f(P), f(P')] \subset \ov T$.  The
smallest such
$C$ is  the
    BBT-constant of $f$, denoted $BBT(f)$.

When $Z$ is a simplicial tree,  we only consider
maps which are linear on each edge.  If $Z$ is a simplicial tree with
a minimal free
action of
$F_k$, every
$F_k$-equivariant  $f:Z\to \ov T$ has bounded backtracking
(see [\BFHg], [\DV], [\GJLL]).

In particular, let $Z_\A$ be the Cayley graph of $F_k$ relative to a free
basis
$\A=\{a_1,\dots,a_k\}$.   The map $f_\A:Z_\A\to \ov T$ sending the
vertex $g$ to
$gQ$ has bounded backtracking, with $BBT(f_\A)\le\sum_{i=1}^k d(Q,a_iQ)$.

If
$w,w'\in F_k$ and $v $ is their longest common initial subword (in the
basis $\A$), then $vQ$ is $BBT(f_\A)$-close to the  segment
$[wQ,w'Q]$ (property
BBT2 of [\GJLL]). This is often used as follows:  if $[Q,wQ]\cap
[Q,w'Q]$ is  long,
then
$vQ$ is  far from
$Q$ and therefore $v$ is  long.

Let $\rho $ be a  ray in $T$.  By [\GJLL, Lemma 3.4],
there is a unique $X=j(\rho )\in\bo F_k$ with the property that
a sequence
$w_n\in F_k$ converges to $X$ if and only if the projection of $w_nQ$ onto
$\rho $ goes off to infinity. The map $j$ is an $F_k$@-equivariant injection
from $\partial  T$ to $\partial   F_k$
satisfying $\da \circ j=j\circ H$. If $\rho $ is an eigenray, then $j(\rho )$
is a fixed point of $\da$. When $\lambda >1$,   every    fixed point of
$\da$ in
$j(\bo T)$ is the image of an    eigenray $\rho $.

The rest of this  section will not  be needed until \S\kern.15em 5.

We suppose   $\lambda >1$. Then
orbits are dense in
$T$ (see [\Pa, Proposition 3.10]), and because   $d(Q,\alpha\mi
(a_i)Q)= \lambda \mi d
(Q,a_iQ)$ there exist bases $\A$ with
$BBT(f_\A)$ arbitrarily small  (this is a special case of [\LLquat,
Corollary 2.3]).

In [\LLquat], we
have associated  a point  $Q(X) \in\ov T\cup\bo T$ to every $X\in \bo
F_k$.  It may be
thought of as    the limit of $g_pQ$ as $g_p\to X$.  Here we shall
mostly be concerned
with whether
$Q(X)$ equals the fixed point $Q$ or not, and we will work with the following
alternative definition of
$Q(X)$.

Given $X\in\bo F_k$, consider the set $B_X$ consisting of  points
$R\in T$ which
belong to the segment
$[Q,wQ]$ for all
$w\in F_k$ close enough to
$X$ (in other  words, $R\in B_X$ if and only if  there exists a neighborhood
$V$ of $X$ in $\ov
F_k$ such that $R\in[Q,wQ]$ for all $w\in V\cap F_k$). It is a
connected subtree
containing no tripod, and there are three possibilities.

If $B_X=\{Q\}$, we define $Q(X)=Q$.

If $B_X$ is unbounded, it is an infinite ray $\rho $ with origin $Q$
and $j(\rho
)=X$.   We define $Q(X)$ as the  point of
$\bo T$ represented   by   $\rho $.

The remaining possibility is that $B_X$ is a closed or half-closed
segment   with
origin
$Q$. We then define
$Q(X)$ as the other endpoint of this segment in $\ov T$ (it may happen that
$Q(X)\notin
B_X$).

It is easy  to check that this definition
of $Q(X)$  coincides with that of [\LLquat]. This implies that  the assignment
$X\mapsto Q(X)$ is
$F_k$-equivariant. Note that in all cases $Q(\da (X))=H(Q(X))$. In particular,
if $Q(X)=Q$, then
$Q(\da ^n(X))=Q$ for every $n\in\Z$.

\snom\lep
\thm{Lemma \stap} Let $T$ be an $\alpha $-invariant \Rt{} as in
Theorem \arbre,   with
$\lambda >1$. Let
$Q\in\ov T$ be the fixed point of $H$. Let
$\A$ be a  basis of $F_k$, and let $f_\A:Z_\A\to \ov T$ be as above
(sending $g$ to
$gQ$). Given
$X\in\bo F_k$, let $X_i$ be its initial segment of length $i$ in the
basis $\A$.
\roster
\item $d(X_iQ, B_X)\le BBT(f_\A)$ for all $i$.
\item If $Q(X)\in\ov T$, then $d(X_iQ, Q(X))\le 2BBT(f_\A)$ for $i$
large enough.
\item If $Q(X)\in\bo T$, the projection of $X_iQ$ onto the ray $B_X$ goes to
infinity as $i\to\infty$.
\endroster
\fthm

\demo{Proof} Suppose $d(X_iQ, Q)>BBT(f_\A)$.
The point
     located on the segment $[Q,X_iQ]$ at distance $BBT(f_\A) $ from
$X_iQ$  belongs to
$[Q,wQ]$ provided $w$ starts with $X_i$, so is  in $B_X$.

Suppose $Q(X)\in\ov T$. Fix  $\varepsilon >0$. For $i$   large, the point
$Q(X)$ is
$\varepsilon $-close to
$[Q,X_iQ]$, and therefore $d(X_iQ, Q(X))\le d(X_iQ, B_X)+\varepsilon $.
Assertion (2) then follows from (1). Assertion (3) is clear.
\cqfd\enddemo

\head {\sect A lemma on asymptotic periodicity} \endhead

Given $\alpha \in\Aut( F_k)$ and $w\in F_k$, we define $w_p=\alpha
^{p-1}(w)\dots\alpha (w)w$ for $p\ge1$ (with $w_1=w$). Note that
$w_r=w_s$   implies
$w_{|r-s|}=1$, and that $w_p=1$ implies
$\alpha ^p(w)=\bigl(\alpha
^{p-1}(w)\dots\alpha (w)\bigr)\mi=w$.

Recall  that $X_n \in \ov F_k$  converges to $X\in\bo F_k$ if and only if $\scp
{X_n}X\to\infty$. The relation  $\scp
{X_n}{Y_n}\to\infty$, between sequences in $\ov F_k$, is transitive
(and does not
depend on the  basis $\A$).

\nom\main
\thm{Lemma \sta} Let $\alpha \in\Aut (F_k)$ and $w\in F_k$. Assume
that for every
$p\ge1$ the elements $w_p$ and $w_p\mi$ are nontrivial and
asymptotically periodic.
Then any
$x\in\ov F_k$ such that $$\lim_{n\to+\infty}\bscp{\ov\alpha ^n(wx)}{\ov\alpha
^{n+1}(x)}=+\infty$$  is  asymptotically periodic.
\fthm

\demo{Proof}
    First suppose that $w$ is not $\alpha $-periodic. Then there exist
$q\ge1$ and $X,Y\in\bo F_k$ such that $\alpha ^{qn}(w)\to X$ and $\alpha
^{qn}(w\mi)\to Y$ as $n\to+\infty$. Write  $\ov\alpha ^{qn}(wx)=\alpha
^{qn}(w)\ov\alpha ^{qn}(x)$. Since $w$ is not periodic, $\alpha
^{qn}(w)$ gets long as $n\to\infty$, and therefore the maximum of
$ \bscp{\alpha ^{qn}(w)}{\ov\alpha ^{qn }(wx)}$ and $\bscp{\alpha
^{qn}(w\mi)}{\ov\alpha ^{qn}(x)}$ goes to infinity with $n$.
This implies that the maximum of
$\bscp{X}{\ov\alpha ^{qn+1}(x)}$ and
$\bscp{Y}{\ov\alpha ^{qn}(x)}$ goes to infinity.
It follows that the limit set
of the sequence
$\ov\alpha ^{qn}(x)$ is contained in $\{\da\mi(X),Y\}$.
Thus   $x$ is  asymptotically
periodic.

Now suppose $\alpha (w)=w$. Then
$\bscp{w\ov\alpha ^n(x)}{\ov\alpha ^{n+1}(x)}$ goes to infinity with
$n$. Fix an
integer $N$. An easy induction shows
$\lim_{n\to+\infty}\bscp{w^N\ov\alpha ^n(x)}{\ov\alpha ^{n+N}(x)}=+\infty$.  As
above, we deduce that for $n$ large at least one of
$\bscp{w^N
}{\ov\alpha ^{n+N}(x)}$,
$\bscp{w^{-N}}{\ov\alpha
^n(x)} $ is large. It follows that
the limit set of
$\ov\alpha ^n(x)$ is contained in $\{w^{-\infty},w^{ \infty}\}$.

The last case is when $w$ is $\alpha $-periodic with period $p\ge2$. Then
$\bscp{w\ov\alpha ^{pn}(x)}{\ov\alpha ^{pn+1}(x)}$ goes to infinity
with $n$, and
our   definition of $w_p$ guarantees
$\bscp{w_p\ov\alpha ^{pn}(x)}{\ov\alpha ^{p(n+1)}(x)}\to\infty$. Since $\alpha
^p(w_p)=w_p$, we reduce to the previous case (replacing $\alpha $ by $\alpha
^p$).
\cqfd
\enddemo

\nom\sousgr
\example{Remark \sta} Suppose $w$ belongs to a finitely generated $\alpha
$@-invariant subgroup $J\inc F_k$. If $x$ is as in Lemma  \main, we
have $\omega
(x)\inc \bo J$ (because all points $X,Y,w^{\pm\infty}$ used in the
proof belong to
$\bo J$).
\endexample

\head  {\sect Limit sets of interior points}\endhead

We first show:

\nom\fevr
\thm{Theorem \sta} Let $\alpha
\in\Aut (F_k)$. If $g\in F_k$ is not $\alpha $@-periodic, the set of limit
points of the sequence
$\alpha ^n(g)$ as $n\to+\infty$ is a periodic orbit of $\da$.
\fthm

\demo{Proof}
We consider an
$\alpha
$@-invariant
\Rt{}
$T$   as  in Theorem  \arbre{}.  Because of Lemma  \fini, we are
free to replace $\alpha $ by a positive   power $\alpha ^r$ whenever
convenient. This
has the effect of replacing $H$ by $H^r$.

The proof is by induction on $k$.
If $g\in\st Q$, with $Q$
a fixed point of $H$, we may use the induction hypothesis (recall
that $\st Q$ is
$\alpha $-invariant and has rank
$<k$).

    We distinguish two cases.

$\bullet$ First suppose $\lambda >1$.
Let
$Q\in\ov T$ be the fixed point of
$H$. If $gQ=Q$, we use induction on
$k$. If
$gQ\neq Q$, and the component
$\C$  of $\ov T\setminus \{Q\}$ containing $gQ$ is fixed by $H$ (in
particular if
$Q\notin T$), that component contains an eigenray  $\rho$ (see
\S\kern.15em 1.c).
Writing
$\alpha ^n(g)Q=\alpha ^n(g)H^nQ=H^ngQ$, we see    that $\alpha ^n(g)$
converges to
$j(\rho )$, a
\fp{}  of
$\da$   (see \S\kern.15em  1.d). Similarly,
$\alpha ^n(g)$ accumulates onto a periodic orbit of $\da$ if the
component $\C $
is
$H$@-periodic.

Assume therefore that $gQ\neq Q$ and $\C $ is not $H$@-periodic. We
will apply Lemma
\main.  Recall that the action of
$\st Q$ on $\pi _0(T\setminus\{Q\})$ has finitely many orbits. Replacing
$\alpha $ by a power, we may assume
that there exists
$w\in\st Q$ with $w \C =H\C $. Define $w_p=\alpha ^{p-1}(w)\dots\alpha
(w)w$ as in   \S\kern.15em 2. The elements $w_p$ are all nontrivial,
because $w_p$
takes
$\C
$ onto  $H^p\C $ (as checked by induction on $p$, using  the equation $\alpha
(w)H=Hw$). Using induction on
$k$, we may assume that $w_p^{\pm1}$ is  asymptotically periodic.

Now we argue as in [\GJLL, p\. 439]. The segments $[Q,wgQ]$ and $[Q, HgQ]$
intersect along a nondegenerate segment. Because $\lambda >1$, the segments
$[H^nQ,H^nwgQ]=[Q,\alpha ^n(wg)Q]$ and $[H^nQ,H^{n+1}gQ]=[Q,\alpha ^{n+1}(g)Q]$
intersect along a segment whose length goes to infinity with $n$. By bounded
backtracking (property BBT2 of [\GJLL],  see \S\kern.15em 1.d), this
implies that the
scalar product
$\bscp{\alpha ^n(wg)}{\alpha ^{n+1}(g)}$ goes to infinity with $n$.  Lemma
\main{} now concludes the proof (with  $\omega (g)\inc \bo\st Q$ by
Remark \sousgr).

$\bullet$  When $\lambda =1$, we view $g$ and $H$ as isometries of
the simplicial
tree $T$. Note that $g$ has at most one fixed point (because  arc stabilizers
are trivial). If
$g$ is hyperbolic (i.e\. it has no
\fp), its translation axis has compact intersection with the  axis of
$H$ (if
$H$ is hyperbolic), or with the set of periodic points of $H$ (if $H$ is
elliptic). Otherwise $g$ would commute with some power  $H^r$ on a
nondegenerate
segment, implying  $\alpha ^r(g)=g$.

First suppose that $H$ is hyperbolic, with axis $A$. Orient $A$ by
the action of
$H$, and consider its two ends $A^-,A^+$.
Choose
$Q\in A$. As
$n$ goes to infinity, the projection of $gH^{-n}Q$ onto $A$ remains far from
$A^-$ (because $g$ does not fix $A^-$). It follows that the
projection of $\alpha
^n(g)Q=H^ngH^{-n}Q$  onto
$A$ goes off to   $A^+$, and $\alpha ^n(g)$ converges to the
    fixed point $j(A^+)$ of $\da$.

If $H$ is elliptic, let $\P$ be the subtree consisting of all $H$@-periodic
points. There exists a    (unique) point $Q$ of $\P$ such that $[Q,gQ]\cap
\P=\{Q\}$ (it is  the point of $\P$ closest to the fixed point of $g$
    if $g$ is elliptic, to the positive end of the axis of $g$ if
$g$ is hyperbolic). Replacing
$\alpha $ by a power, we may assume
$HQ=Q$.

If
$gQ=Q$, we use induction on $k$. If not, we let $\C $ be the component
    of $T\setminus \{Q\}$ containing $gQ$, and we find $w\in\st Q$ with
$w\C =H\C $,
as in the
     proof when $\lambda >1$. The components  $\C_n=H^n\C=w_n\C $ are all
distinct, because $Q$ was chosen  so that the germ of $[Q,gQ]$ at $Q$ is not
$H$@-periodic.
Since both   points
$\alpha ^n(wg)Q=H^n wgQ$ and $\alpha ^{n+1}(g)Q=H^{n+1}gQ$  belong to
$\C_{n+1}$, the scalar product
$\bscp{\alpha
^n(wg)}{\alpha ^{n+1}(g)}$ goes to infinity and Lemma \main{} applies.
\cqfd
\enddemo

The argument that concludes this proof will be used again. We may  state
it as follows.

\nom\maine
\thm{Lemma \sta} Let $T$ be a minimal simplicial $F_k$-tree with trivial edge
stabilizers. Given
$Q\in T$ and distinct components $\C_n$ of  $T\setminus\{Q\}$, there exists a
sequence of  numbers
$m_n\to\infty$ such that, if $a_n,b_n$ are elements of $F_k$ with $a_nQ$
and $b_nQ$ both belonging to $\C_n$, then $\scp{a_n}{b_n}\ge m_n$. \cqfd
\fthm

   From Theorem \fevr{} we deduce:

\nom\trick
\thm{Corollary \sta}   Let $\alpha \in\Aut(F_k)$ and
       $w\in F_k$. If the sequence $w_n=\alpha ^{n-1}(w)\dots\alpha
(w)w$ is not periodic (in particular   if $w$ is not $\alpha
$-periodic), its limit
set  is a periodic orbit of
$\da
$.
\fthm

\demo{Proof} Extend
$\alpha
$ to
$\beta
\in\Aut(F_k*\Z)$ by sending  a generator $t$ of $\Z$ to $wt$. Then $\beta
^n(t)=w_nt$. The map $n\mapsto w_n$ is injective, because otherwise
$w_n$ would be
periodic. Thus the length of $w_n$ goes to infinity, implying that
the sequences
$w_n$  and
$\beta ^n(t)$ have the same limit set in $\bo(F_k*\Z)$. That set is
contained in
$\bo F_k$, and is a periodic orbit of $\db$ by Theorem \fevr{}. It is
therefore a
periodic orbit of $\da$. \cqfd\enddemo

The following observation will be useful in \S\kern.15em 9.

\nom\bij
\thm{Proposition \sta} Let $T$ be an $\alpha $-invariant \Rt{}. Suppose that
$\lambda >1$ and   the fixed point $Q\in\ov T$ of $H$ has trivial
stabilizer. Then
$\alpha $ has no nontrivial periodic element in $ F_k$. There is a bijection
$\tau $ from $\pi _0(\ov T\setminus \{Q\})$ to the set of attracting
periodic points
of $\da$, with $\tau \circ H=\da\circ\tau $.
\fthm

\demo{Proof} If $\alpha ^q(g)=g$, then $H^q gQ=\alpha ^q(g)H^qQ=gQ$,
so $g\in\st
Q=\{1\}$. This proves the first assertion.
By [\GJLL], $\da$ has finitely many
periodic points, each of them attracting or repelling.

Since $\st Q$ is trivial, $\ov T\setminus \{Q\}$ has finitely many
components, each of them $H$-periodic. As seen in the proof of Theorem
\fevr, the asymptotic behavior of
a sequence $\alpha ^n(g)$ depends only on the component $\C$
containing $gQ$.  The
attracting periodic point of $\da $ associated to
$\C$   is
$\tau (\C)=j(\rho )$, where $\rho $ is the $H$-periodic ray contained in $\C$.
\cqfd\enddemo

\example{Remark \sta} Suppose $\lambda >1$, but   don't  assume that $\st
Q$ is trivial. The
proof of  Theorem
\fevr{} shows that either $\omega (g)\inc \ov{\st Q}$, or $\alpha
^n(g)$ accumulates onto
a periodic orbit of $\da$ associated to $H$-periodic rays.
\endexample

\head \sect The simplicial case\endhead

Recall  that $\alpha \in\Aut(F_k)$ has asymptotically periodic dynamics if
every $X\in\ov F_k$ is asymptotically  periodic.  This section is devoted
to the proof of the following result.

\nom\simpli
\thm{Theorem \sta}  Let $\alpha \in\Aut(F_k)$. Assume  that all
restrictions  of $\alpha $ to $\alpha $-invariant subgroups of rank
$<k$ have asymptotically periodic dynamics.   If there exists a simplicial
$\alpha
$-invariant tree
$T$  as  in Theorem  \arbre{}, then $\alpha $ has asymptotically
periodic dynamics.
\fthm

\nom\coroo
\thm{Corollary \sta} Polynomially growing  automorphisms of $F_k$ have
asymptotically periodic dynamics.
\fthm

\demo{Proof of Corollary \coroo} By induction on
$k$. Given $\alpha
\in\Aut(F_k)$, consider  $T$   as in Theorem \arbre{}, with length
function $\ell$.
Choose $w\in F_k$ with
$\ell(w)>0$. If
$\alpha
$ is polynomially growing,   then
$\ell(\alpha ^n(w))$   has subexponential growth (as  translation length is
bounded from above  by a constant times   word length). Since
$\ell(\alpha ^n(w))=\lambda
^n\ell(w)$, we have  $\lambda =1$. By Assertion (c) of Theorem \arbre, the
tree $T$ is simplicial. Theorem
\simpli{} applies  since any restriction of $\alpha $ is polynomially
growing.
\cqfd\enddemo

    To prove Theorem \simpli, we argue as in the
proof of Theorem \fevr{} when $\lambda =1$ (keeping the same notations). If the
isometry
$H$ is hyperbolic, then
$\da ^n(X)$  converges to the
    fixed point $j(A^+)$ for every $X\neq j(A^-)$, so we assume that $H$
is elliptic.

First suppose that
$X $  belongs to
$\bo\st R$ for a (unique) point $R\in T$.  Let $Q$ be
the  point of $\P$ closest to $R$. Replacing
$\alpha $ by a power, we may assume $HQ=Q$.
If $R=Q$ (i.e\. if $R$ is fixed by $H$), we use induction on $k$. If
$R\neq Q$, we let  $\C$ be the component of
$T\setminus \{Q\}$ containing $R$ and we choose $w\in\st Q$ with $w\C=H\C$
(replacing
$\alpha $ by a power if needed).

Let $g_p$ be a sequence in $\st R$ converging to
$X$. Since $g_pQ\in \C$, we can apply Lemma \maine{} for each value
of $p$, with
$\C_n=H^n\C$,
$a_n=\alpha ^n(wg_p)$, $b_n=\alpha ^{n+1}(g_p)$. We find $\bscp{\alpha
^n(wg_p)}{\alpha ^{n+1}(g_p)}\ge m_n$, and therefore
$\bscp{\da^n(wX)}{\da^{n+1}(X)}$ tends to infinity as $n\to+\infty$.
    Lemma
\main{} concludes the proof, since $w$ belongs to  an $\alpha $-invariant
subgroup of rank $<k$.

If $X$ does not belong to any $\bo\st R$, then, since  $T$ is simplicial and
arc stabilizers are trivial,
$X=j(e)$ for some end
$e$ of
$T$.  We may assume that
$e$ is not $H$@-periodic. We apply Lemmas \maine{} and \main{}
as above, using the   point
$Q\in \P$ closest to $e$ and the component $\C$ of $T\setminus \{Q\}$
containing
$e$. The existence of this point $Q$ is not obvious, however, because
$e$ could be a limit of periodic points of $H$ whose period tends to infinity.
But this is ruled out by the following fact, and the proof  of Theorem
\simpli{} is complete.

\nom\finiper
\thm{Lemma \sta}
Let $\alpha ,T,H$ be as in  Theorem \arbre{}. Assume $\lambda =1$.
There exists $M$ (depending only on $k$) such that all
periodic points of $H$ have period less than $M$.
\fthm

\demo{Proof} If the action of $F_k$ on $T$ is free,  some fixed power of
$H$ is hyperbolic or is the identity (because $H$ lifts an automorphism of the
finite graph $T/F_k$). We assume from now on that the action is not free.

Suppose that a vertex $Q$, with $\st Q$ nontrivial, is $H$@-periodic with
period $q$. Since $\st Q$ is invariant under $\alpha ^q$, we know
that $\bo\st Q$
contains a periodic point $X$ of $\da^q$. Viewed as a periodic point of $\da$,
the point $X$ has period divisible by $q$ (because stabilizers of distinct
vertices $H^iQ$, $1\le i\le q$, have disjoint boundaries), and this forces
$q\le M_k$ by  [\LLdeux] (see  Theorem \rappel).

We have now obtained a bound for periods (under $H$) of vertices with
nontrivial
stabilizer. Since the action of $F_k$
on $T$ is minimal and  not free, any vertex $Q$ with trivial stabilizer belongs
to a segment $[Q_1,Q_2]$, with $\st Q_1$ and $\st Q_2$ nontrivial but $\st R$
trivial for every interior point $R$ of $[Q_1,Q_2]$. If $Q$ is
$H$@-periodic, so
are $Q_1$ and $Q_2$ because vertices with trivial stabilizer have
finite valence in
$T$. We  know how to bound the period of $Q_1$ and $Q_2$, and this  leads to
a bound for the period of $Q$.
\cqfd\enddemo

\head\sect A reduction  \endhead

   From now on we will consider an $\alpha $-invariant \Rt{} $T$  with
$\lambda >1$,
as in Theorem \arbre. We always denote by
$Q$ the fixed point of $H$ (in $\ov T$). We fix $X\in\bo F_k$, and we study the
sequence $\da^n(X)$.

For the reader's convenience, we now give a quick overview of \S\S
\kern.15em 5-8.

In \S\kern.15em 3,  we studied the sequence $\alpha ^n(g)$ by applying
powers of
$H$ to the point $gQ$. Instead of  $gQ$, we now consider the point
$Q(X)$ introduced
in [\LLquat] (see \S\kern.15em 1.d). It is either a point of  the
metric completion  $\ov
T$, or an end of
$T$. It may be thought of as    the limit of $g_pQ$ as $g_p\to X$.

When $Q(X)\neq Q$, a naive argument works: the behavior of
$\da^n(X)$ is the same as that of $\alpha ^n(g)$ for $g\in F_k$ close
to $X$.  If $X\in
\bo\st Q$, we use induction on $k$. The hard case is when
$Q(X)=Q$ but $X$ does not belong to $\bo\st Q$. This happens when $X$ is a
repelling fixed point of $\da$,
and the task here is indeed to show that this is essentially the only
possibility
(Theorem 5.1).

The difficulty  lies in the fact that
cancellation within the infinite word $X$ due to the   contracting
nature of the
combinatorics of a repelling $X$ can be mixed with cancellation in the
initial subword of $X$ coming from accumulations of finite elements
from $\st Q$.

To overcome this  difficulty, we 
consider an improved  relative   train track map $f_0:G\to G$
representing $\alpha
$ (in the sense of [\BFHt]), with an exponential top stratum, and we
use  the   invariant tree $T$   associated to  that stratum as in [\GJLL].

If
there is an indivisible Nielsen path $\eta $
meeting the top stratum, then for each lift
$[a,b]$ of
$\eta $ in the universal covering $\wtilde G$ we create a new edge
between $a$ and $b$
(this is similar to the addition of Nielsen faces  in [\Luconj]). The
resulting space
$\tg$ is a cocompact
$F_k$-space, but we also consider the non-proper space $\tg\PF$
obtained from $\tg$ by
contracting all edges not in the top stratum (it resembles  the
coned-off Cayley
graph of [\Far]). We show that, when
$Q(X)=Q$, the point at infinity $X$ behaves like a repelling point
for the train track
map acting on $\tg\PF$ (Lemma 8.2).

We also  need to consider a train track map $f'_0:G'\to G'$ representing
$\alpha \mi$, which is
paired with $f_0$ in the sense of [\BFHt, \S\kern.15em 3.2]. Because
of the pairing, the
spaces
$\tg\PF$ and $\tg'\PF$ are   quasi-isometric. It follows that $X$
behaves like an
attracting point in $\tg'\PF$, and there is a legal ray going out
towards $X$ in
$\tg'$ (Lemma 8.4). We then conclude by applying these facts to all points
$\da^{-n}(X)$.

We now proceed to reduce Theorem \pt{} to the following result, whose proof
will be completed in \S \kern.15em 8.

\nom\leredu
\thm{Theorem \sta}
Let $\alpha \in\Aut(F_k)$. If there is no simplicial $\alpha
$-invariant tree
(with trivial arc
stabilizers),
there is an $\alpha $-invariant \Rt{} $T$  with $\lambda >1$ such
that, for every
$X\in\bo F_k$, at least one of the following holds:
\roster
\item $X\in\bo\st Q$.
\item $Q(X)\neq Q$.
\item There exist $q\ge1$ and  $w\in\st Q$ such that $X$ is a fixed point of
$w\da^q=\bo( i_w\circ\alpha^q)
$.
\endroster
\fthm

      See \S \kern.15em 1.d for the definition of  $Q(X)
\in\ov T\cup\bo T$, and recall that $i_w(g)=wgw\mi$.

Assuming this theorem, we prove  Theorem \pt{}
by induction on
$k$. Consider $\alpha \in\Aut(F_k)$. If there is a simplicial $\alpha
$-invariant
tree, we apply Theorem
\simpli.  If not, we fix $X\in\bo F_k$ and we  consider the three
possibilities of
Theorem
\leredu.

$\bullet$ If $X\in\bo\st Q$, we use the induction
hypothesis (recall that $\st Q$ is $\alpha $-invariant, with  rank $<k$).

$\bullet$ If $Q(X)\neq Q$,
let $\C$ be the component
of $\ov T\setminus \{Q\}$ containing $B_X\setminus \{Q\}$.
    Choose a basis $\A$ of $F_k$ with $BBT(f_\A)$ small with
    respect to $d(Q(X), Q)$. By
Lemma
\lep,  we have $X_iQ\in\C$ for $i$ large.
    Replacing
$\alpha $ by a power if necessary, we may assume (as done in the
previous sections) that
there exists
$w\in\st Q$ such that $H\C=w\C$. We argue as in the proof of Theorem
\fevr, ``uniformly''. There are two cases.

If $w=1$, let
$\rho
$ be the  eigenray of
$H$ contained in $\C$. The element $j(\rho )\in\partial F_k$ is a
\fp{} of $\da$,
and we   show
$\lim_{n\to+\infty}\da^n(X)=j(\rho )$.

The  distance from $Q$ to the projection of $X_iQ$ onto $\rho $ is
    greater than some $\varepsilon   >0$ for $i$ large.
      It follows that the distance from $Q$ to the projection of
$\alpha ^n(X_i)Q=H^n(X_iQ)$ onto $\rho $ is greater than $\lambda
^n\varepsilon  $,
     independently of  $i$.
By bounded backtracking, as in the proof of Lemma 3.4
in [\GJLL], we obtain that $\bscp{{j(\rho )}}{\alpha ^n(X_i)}$ goes
to infinity as
$n\to+\infty$,  uniformly in $i$.   This
shows
$\lim_{n\to+\infty}\da^n(X)=j(\rho )$, as required.

If $w\neq1$, we  observe that the length of $[Q,wX_iQ]\cap[Q,HX_iQ]$ is bounded
away from $0$, so that
$\lim_{n\to+\infty}\bscp{\ov\alpha ^n(wX_i)}{\ov\alpha
^{n+1}(X_i)}=+\infty$   uniformly in $i$. Therefore $\bscp{\ov\alpha
^n(wX)}{\ov\alpha ^{n+1}(X)}$ goes to infinity and Lemma \main{}  applies.

$\bullet$ Replacing $\alpha $ by $\alpha ^q$, we may assume   that $X$
is a fixed point of
$\bo \beta
$, with $\beta =i_w\circ\alpha$.
Write $\bo\alpha ^n(X)= w^-_nX$ with   $w^-_n=\alpha
^{n-1}(w\mi) \dots\alpha  (w\mi) w\mi\in \st Q$. If there exist $r\neq s$ with
$w^-_r=w^-_s$, then $X$ is $\da
$-periodic and we are done.
We therefore assume $|w^-_n|\to\infty$. We may also assume that   the
cancellation
in the product
$w^-_nX$ is bounded, since otherwise $X\in\bo\st Q$. It follows that the
sequence
$\da^n (X)$ has the same limit points as the sequence $w^-_n$, and we
conclude by
Corollary \trick.

\head\sect The train track and the tree\endhead

In this section we let the map
    $f_0:G\to G$ be an improved relative train track representative for
(a power of) $\alpha
$, in the sense of [\BFHt, 5.1.5]. Basic references are [\BH] and
\S \S\kern.15em 2.5, 5.1 of [\BFHt].

We assume that the top stratum $G_t$ is exponentially growing.  Edges in
 $G_t$ are called {\it top edges\/}, edges in the other   strata
are called {\it zero edges\/}.
An illegal turn in $G$ is
a turn between two  top edges     where folding occurs when some
power of
$f_0$ is applied (see [\BH]).

There is at most one indivisible Nielsen path (INP) $\eta  $ in $G$
meeting the top
stratum (see [\BFHt, 5.1.5 and 5.1.7] for its properties). It has
precisely one illegal
turn (its tip), and $f_0(\eta )$ is homotopic to $\eta
$ relative to its endpoints (which are fixed by $f_0)$. If $\eta $ is
a loop, there
is no cancellation in $\eta ^2$ and the  turn at the midpoint of
$\eta ^2$ is legal.

Our actual object of study will be   a
 {\it train track with
shortcut\/} $\G$, obtained from $G$ by adding a new zero edge if needed. This
is  a very special case of
the ``train tracks with Nielsen faces'' introduced in [\Luconj].

\subhead  \subsect  Creating a shortcut
\endsubhead

    If there is no INP $\eta $ as above,
we let $\G=G$. If there is  one, we first subdivide $G$  (if needed)
to ensure that
the endpoints of $\eta $ are vertices. Then we enlarge
$G$ into
$\G$ by adding a new zero edge
$e$ with the same endpoints as
$\eta
$ (these endpoints are equal if the stratum is geometric).
    The   map $f_0$, originally defined  on $G$,  extends to $\G$ (it is
defined as  the identity on
$e$).

We fix a retraction
$r_0:\G\to G$ mapping $e$ onto $\eta $ in  a locally injective way.
It induces a
homomorphism   $(r_0)_*:\pi _1(\G)\to \pi _1(G) \simeq F_k$, and we
let $\wtilde\G$
be the corresponding
$F_k$-covering (it consists of the universal covering $\wtilde G$ of
$G$, together
with lifts of
$e$; if a $2$-cell is attached to $\G$ along $\eta \cup e$, then $\tg$ becomes
the $1$-skeleton of the universal covering).

We define  top   and zero edges on $\tg$ as on
$\G$, we lift $r_0$ to a retraction $r$, and
we define $f$ to be the lift of $f_0$ to $\wtilde \G$ that  satisfies $\alpha
(w)f=fw$ for every deck transformation $w\in F_k$.

\subhead  \subsect   Legal paths
\endsubhead

We  define an   {\it illegal turn  in $\G$\/}
as either an illegal turn in $G$, or
     a turn between $e$ and a
    top edge   which becomes degenerate or illegal in
$G$ when
$e$ is replaced by $\eta $ (recall   that $\eta $ starts and
ends with segments contained in the top stratum $G_t$; if an endpoint
$x$ of $\eta $
lies in the interior of a top edge of $G$, exactly one of the turns
at the corresponding
end of
$e$  is illegal). A turn in $\tg$ is illegal if it projects onto an
illegal turn.

Paths   will always be assumed to start and end at vertices. A {\it
legal path\/}, in
$
\G$ or in
$\tg$, is a path with no illegal turn. Note that    $r_0$ maps a
legal path  locally
injectively,   and the image has  illegal turns   only at the tip of
$\eta $. Given
vertices $x,y\in\tg$, we let $ILT(x,y)$ be the minimum number of
illegal turns on
paths from
$x$ to $y$ in  $\tg$.

Now consider the action of $f_0$ on
a legal path $\gamma \inc\G$. Write $\gamma $ as a concatenation of
copies of $e$
and legal paths in $G$. When $f_0$ is applied, the paths in $G$ remain
legal and nontrivial (after reduction), by property (RTT-ii) of
[\BH], and the turns
at the endpoints of $e$ remain legal. It follows that $f_0$ and $f$ map  legal
paths to  legal paths. In particular, we have  $ILT(f(x),f(y))\le ILT(x,y)$.

\snom\iterleg
\thm{Lemma \stap{} [\BFHt, \Luconj]}  Let $x,y$ be vertices of $\tg$. For $p$
large enough, there is a legal path joining $f^p(x)$ to
$f^p(y)$. More precisely, given $\delta >0$, there exists $p$ such that
$ILT(f^p(x),f^p(y))\le \delta ILT(x,y)$ for all vertices $x,y\in\tg$.
\fthm

\demo{Proof}
By Lemma 3.2 of [\Luconj], or Lemma
4.2.6 of [\BFHt], 
 there exists $p$ such that $f^p(x)$ and
$f^p(y)$  may be joined  in $\wtilde G$  by a    legal concatenation of
legal paths and INP's.  This implies the first assertion  
(note that the definition of a legal path used here is slightly
stronger than the one
used in [\Luconj], but it is easy to make the path legal in the
stronger sense by increasing   $p$).

By Remark 3.3 of [\Luconj], there exists
$p$ such that
$ILT(f^p(x),f^p(y))=0$ if
$ ILT(x,y)=1$. This easily implies $ILT(f^p(x),f^p(y))\le \frac12\,
{ ILT(x,y)} $ for
all  $x,y $. The extension to an arbitrary $\delta $ is immediate.
\cqfd\enddemo

\subhead  \subsect  The PF-metric and the invariant tree
\endsubhead

Since the top stratum $G_t$ is assumed to be exponentially growing,
the transition matrix
associated to $G_t$ has a Perron-Frobenius eigenvalue $\lambda >1$. Using
components of an eigenvector, one assigns a PF-length $|E|\PF$ to
each top edge $E\inc
\G$, with the property that the total PF-length of $f_0(E)$ is
$\lambda|E|\PF$. More
generally, $f_0$ multiplies the PF-length of legal paths by $\lambda $.

On each of the spaces $\wtilde G$ and $\tg$, we  define  a {\it
PF-pseudo-distance}
$d\PF$  by giving top edges their PF-length and assigning  length $0$
to the zero
edges (as usual, the distance between two points is the length of the
shortest path
joining them). Note that these two distances differ on
$\wtilde G$: the endpoints of a lift of
$\eta $ have PF-distance $0$ in $\tg$ but not in $\wtilde G$.  The
associated metric
spaces
$\wtilde G\PF$ and $\tg\PF$ are geodesic (because there is a positive
lower bound for
the PF-length of top edges), but in general   not proper.

We also define the {\it simplicial distance\/} $d$, for which top edges have
their PF-length and zero edges have length $1$. This makes $\wtilde G$
and $\tg$ into
    proper geodesic spaces, with   cocompact actions of $F_k$. Note
that  $d\PF\le
     d$.

Now recall the construction of an $\alpha $-invariant \Rt{} $T$
with  trivial arc
stabilizers given in [\GJLL].
The pseudo-distance $d_{PF}$ on the tree $\wtilde G$    satisfies
$d\PF(f(x),f(y))\le\lambda d\PF(x,y)$, with equality if $x$ and $y$
are joined by a
legal path.   We define
$$d_\infty(x,y)=\lim_{p\to\infty}\lambda ^{-p}d\PF(f^p(x),f^p(y))$$
and we let $T$
be the metric space associated to this pseudo-distance. We denote by
$\pi $ the quotient
    map $\pi :\wtilde G\to T$.
    The map $f$ induces a
homothety  $H:T\to T$, with stretching factor $\lambda $, such that  $\pi
\circ f=H\circ
\pi $. 

It is proved in [\GJLL] that $T$
 satisfies the conditions of Theorem
\arbre{}. We say that $T$ is the {\it invariant tree associated to the train
track map $f_0:G\to G$\/} (and its exponential top stratum; when the top
stratum of
$G$ is polynomially growing, one obtains a simplicial invariant tree simply
by collapsing all zero edges of
$\wtilde G$).

    Note that   the endpoints of a given lift of $\eta
$ have the same image in $T$, since their images by $f$ also bound a
lift of $\eta
$. This implies that,  if  we apply  the
above construction  to $d\PF$ on $\tg $ rather than on $\wtilde G $,
we get the same space $T$, and the quotient map $i:\tg\to T$ extends $\pi $.

    The
maps
$i:\tg\to T$ and
$\pi :\wtilde G\to T$  are
$F_k$-equivariant, and $1$-Lipschitz for the respective $d\PF$ (hence
also for the
simplicial metrics). Also note that  $i$ and $\pi \circ r$ are uniformly close.

The map $\pi $, defined on the tree $\wtilde G $, sends legal paths
PF-isometrically into
$T$, and (see
\S\kern.15em 1.d) it    has bounded backtracking (the image of a
geodesic segment $[x,y]$
is close to $[\pi (x),\pi (y)]$). We show that
$i:\wtilde\G_{PF}\to T$ has similar properties.  This is a quite
different  statement,
because geodesics in
$\wtilde\G_{PF}$ can be very different from simplicial geodesics.

\snom\legalge
\thm{Lemma \stap{} [\Luconj, Lemma 3.9 and Proposition 4.2]} A legal
path $\gamma
\inc\tg$ is PF-geodesic and maps to
$T$ PF-isometrically. The map $i:\wtilde\G_{PF}\to T$ has bounded backtracking.
\fthm

\demo{Proof} Let $x,y\in\wtilde G$ be the
endpoints of   $\gamma $.
Let $L$ be the PF-length of $\gamma $, and
$m$ the number of lifts of
$e$ contained in $\gamma $. Then $f^p(\gamma )$ is a legal path
     with PF-length $\lambda ^pL$,
containing   only
$m$ lifts of $e$. Its image by $r$ is embedded in $\wtilde G$ and has  PF-length
$\lambda ^pL+m|\eta |\PF$. Letting $p\to\infty$ shows that  $i (x)$ and
$i (y)$ have distance
$L$ in
$T$. Since
$i$ is 1-Lipschitz for the PF-metric, we see that $\gamma $ is
PF-geodesic and $i_{|\gamma }$ is  isometric for the PF-metric.

To prove bounded backtracking, it suffices to
     consider a PF-geodesic $\gamma
=[x,y]$ in
$\wtilde\G$ with $i(x)=i(y)$, and  to bound the diameter of $i(\gamma )$. We
may assume $x,y\in\wtilde G$. Let $\delta $ be the geodesic $[x,y]$ in the tree
$\wtilde G$. We know that $\pi $ has bounded backtracking, so
     $\pi (\delta )$ is close to $i(x)$.
Consider $z\in\gamma $
with $i(z)$ far from $i(x)$. Since $i$ and $\pi \circ r$ are close,
we have  $r(z)\notin\delta $. Because $\wtilde G$ is a  tree,
     $z$ belongs to a subarc $[x',y']\inc\gamma $ with $r(x')=r(y')\in\delta $.
The
points $x'$ and $y'$ are close in $\tg$ for the simplicial distance,
hence also for
the PF-distance. Since $\gamma $ is PF-geodesic (and $i$ is
$1$-Lipschitz  for the
PF-metric),
$i(z)$ is close to $i(x')$, hence to $\pi (r(x'))$, hence to $i(x)$.
\cqfd\enddemo

\subhead \subsect  Elliptic elements \endsubhead

A loop in $G$  represents a conjugacy class in $\pi
_1(G)\simeq F_k$. So does a loop in $\G$, through the   homomorphism
$(r_0)_*:\pi _1(\G)\to \pi
_1(G)\simeq F_k$.

\snom\ellipt
\thm{Lemma \stap} Given a conjugacy class $[u]$ in $F_k$, the following are
equivalent:
\roster
\item
$[u]$ may be represented by a loop in the zero part of $\G$;
\item
$[u]$ may be represented by a loop in the zero part of $G$, or $[u]$ is $\alpha
$-invariant (i.e\. $\alpha (u)$ is conjugate to $u$);
\item
$[u]$ is elliptic in the associated tree  $T$ (i.e\. $u$ fixes a point in
$T$).
\endroster
\fthm

\demo{Proof}
For $p>0$,    each of the three conditions   holds for $[u]$ if and
only if it holds for     $\alpha ^p([u])$ (for the first two
conditions, use Scott's
lemma [\BFHt,  6.0.6] to argue that
$[u]$ may be represented by a loop in the zero part if $\alpha
^p([u])$ does). By
Lemma
\iterleg{} we may assume that
$[u]$ is represented by a legal loop $\gamma $ in $\G$. By Lemma \legalge, the
translation length $\ell(u)$ of
$[u]$ in $T$ is  then the PF-length of $\gamma $.

If this length is positive, then $u$ is hyperbolic in $T$. It      cannot be
represented by a loop in a zero part, or be $\alpha $-invariant
(since $\ell(\alpha
(u))=\lambda
\ell(u)$).
     If the length is $0$, then   $[u]$ is elliptic and $\gamma $ is
contained in the
zero part of
$\G$. If $[u]$   cannot be represented by a loop in the zero part of $G$, then
by [\BFHt, 6.0.2]
$u$ is (conjugate to) a power of $\eta $ (which is a loop). It
follows that $[u]$ is $\alpha
$-invariant.
\cqfd\enddemo

We may be more specific, working with elements of $F_k$ rather than
with conjugacy
classes. Let $Z_i$ be the   non-contractible
components of the zero part of
$G$. Images of $\pi _1(Z_i)$ in
$F_k$ are free factors. By [\BFHt, 5.1.5], these factors are $\alpha
$-invariant (up to conjugacy). Not much changes when we pass from $G$
to $\G$. If
$G_t$ is not geometric, zero components of $\G$ yield the same
subgroups [\BFHt,
5.1.7]. In the geometric case, there is one more component, whose
image in $F_k$ is
the cyclic group generated by the loop $\eta $ (which runs once
around the closed INP).

If $Z$ is a non-contractible zero component of $\G$, and $\hat Z$ is
a component of
its preimage in $\tg$, the image of $\hat Z$ in $T$ is a point
$i(\hat Z)$. By Lemma
\ellipt{} the image of $\pi _1(Z)$ in $F_k$ is the full stabilizer of
$i(\hat Z)$,
and every non-trivial stabilizer arises in this way: non-trivial stabilizers of
points of $T$ are precisely conjugates of fundamental groups of
non-contractible
zero components of
$\G$.

\head \sect Geometry on the train track \endhead

 As in the previous section, $f_0:G\to G$ is a train track map, $\G$ is
the train track with shortcut (obtained by adding   an edge if there is an
  INP), and
$T$ is the associated invariant
\Rt.

\subhead \subsect  Spaces quasi-isometric to trees \endsubhead

If $c$ is a continuous map from $[0,1]$ to a tree, the image of $c$
contains the
geodesic segment between $c(0)$ and $c(1)$. Spaces quasi-isometric to
a tree have a
similar property.

\snom\truc
\thm{Lemma \stap} Let $Y$ be a geodesic metric space quasi-isometric to a tree.
There exists $C$ such that, if $\gamma $ is a geodesic segment and
$\gamma '$ is
any path with the same endpoints, then every $P\in \gamma $ is
$C$-close to some
$Q\in\gamma '$.
\fthm

\demo{Proof} Let $f:Y\to S$ be a quasi-isometry to a tree. Since
$f(\gamma )$ is a
quasi-geodesic,
$f(P)$ is close to some $R$ on the geodesic $\gamma _0$ joining the
$f$-images of the
endpoints of
$\gamma $, and $R$ is close to the image of  some $Q\in \gamma
'$. The distance from $P$ to this point $Q$ is bounded.
\cqfd\enddemo

This basic fact may be extended in several ways.

If the endpoints of $\gamma '$ are only $C$-close to those of $\gamma $, then
every $P\in\gamma $ at distance at least $2C$ from the endpoints is
$C$-close to
$\gamma '$. If   points $P_1,P_2,\dots$ appear in this order on
$\gamma $, each  at
distance at least
$2C$ from the previous one, we may assume that the associated points $Q_i $
    appear in the same order on
$\gamma '$ (construct $Q_i$
inductively).

Lemma \truc{} also holds if $\gamma ,\gamma '$ are (possibly infinite)
quasi-geodesics with the same endpoints, with $C$ depending only on the
quasi-geodesy constants of $\gamma $.

\subhead \subsect  $K$-PF-geodesics  \endsubhead

Let $\tg$ be as above. It is equipped with the simplicial distance $d$, the
PF-distance $d\PF$, and the distance $d_\infty$.  Recall  that
$d_\infty\le d\PF\le
d$. The space $(\tg, d)$ is proper and quasi-isometric to $F_k$,
hence to a tree. We
identify $\bo\tg$ with $\bo F_k$. We fix $C$ as in Lemma \truc.

The words ``geodesic'' and ``quasi-geodesic'' will always refer to
the simplicial
metric. We write
$K$-quasi-geodesic instead of
$(K,K)$-quasi-geodesic. A  geodesic relative to  $d\PF$ will be
called a PF-geodesic. If a PF-geodesic is also  $K$-quasi-geodesic
(with respect to $d$), we call it   a $K$-PF-geodesic. Recall that
legal paths are PF-geodesics.

\snom\geodspec
\thm {Lemma \stap } If $K$ is large enough, the following  properties hold:
\roster
\item There exists a
K-PF-geodesic between any two points of $\tg\cup\bo\tg$.
\item If there exists
    a legal
path between two points of $\tg$, there exists a legal
K-quasi-geodesic between them.
\endroster
\fthm

\demo{Proof}
Given a  geodesic segment $\gamma $, first choose a PF-geodesic $\gamma '$ with
the same endpoints. Place points
$P_i,Q_i$ as above (after Lemma \truc), with $d(P_i,P_{i+1})= 2C$ (we assume
that the points
$Q_i$ are vertices). Replace the
segment of
$\gamma '$ between
$Q_i$ and
$Q_{i+1}$ by another PF-geodesic segment,
with   simplicial length as small as possible.
The resulting curve   is a PF-geodesic,  and   it is uniformly quasi-geodesic
because the set of pairs
$(Q_i,Q_{i+1})$   is finite up to deck transformations ($\tg$ is a
locally finite  graph).  If
$\gamma $ is  infinite, we apply the usual diagonal argument.

If there is a legal $\gamma '$, we choose the replacement of  $[Q_i, Q_{i+1}]$
among legal paths with the same initial and final edges (so that the
turns at the
$Q_i$'s remain legal).
\cqfd\enddemo

\snom\llq
\thm{Lemma \stap} Let $\gamma\inc \tg$ be a   quasi-geodesic ray, with point at
infinity
$X$. If $Q(X)\in\ov T$, then $Q(X)$ belongs to the closure of $i(\gamma )$. If
$Q(X)\in\bo T$, then $i(\gamma )$ contains a ray going out to $Q(X)$.
\fthm

\demo{Proof}  Fix $\varepsilon >0$. Choose a basis $\A$ such that
$f_\A:Z_\A\to T$ has
backtracking less than
$\varepsilon $ (see  \S \kern.15em 1.d). As in [\LLquat, proof of
3.1], subdivide
$\tg$ to get edges all of simplicial length $<\varepsilon $  and construct an
equivariant map $\zeta :\tg\to Z_\A$ such that  $f_\A\circ \zeta $ is
$2\varepsilon
$-close to
$i$. As $Z_\A$ is a  tree, the image of $\gamma $ in $Z_\A$ contains
a ray $\rho $
going out to
$X$.  If $Q(X)\in\ov T$, it is   $2\varepsilon $-close to
$\ov{f_\A(\rho )}$ by Lemma
\lep, hence
$4\varepsilon
$-close to
$i(\gamma )$. If $Q(X)\in\bo T$, then $f_\A(\rho )$ contains a ray going out to
$Q(X)$ by Lemma
\lep.
\cqfd\enddemo

Let $BBT(i)$ denote the backtracking constant of $i :\wtilde\G_{PF}\to
T$ (see Lemma
\legalge).

\snom\coror
\thm{Corollary \stap}
If   $\gamma \inc\tg$ is a K-PF-geodesic ray, with point at infinity
$X$, and   $i$ maps the origin of $
\gamma
$  to $Q(X)\in\ov T$,  then
$i(\gamma )$ is contained in the $BBT(i)$-ball centered at $Q(X)$.
\fthm

\demo{Proof} Suppose a point $x\in\gamma $ is mapped by $i$ at
distance $>BBT(i)$
from $Q(X)$. Since $\gamma $ is PF-geodesic, the point at distance
$BBT(i)$ from
$i(x)$ on the segment
$[Q(X),i(x)]$ separates $Q(X)$ from $i(y)$ for $y\in\gamma $ closer to $X$ than
$x$. This contradicts Lemma \llq.
\cqfd\enddemo

\subhead \subsect An inequality   \endsubhead

Recall that  $ILT(x,y)$ is  the minimum number of illegal turns in any
    path from $x$ to $y$ in  $\tg$.

\snom\borneilt
\thm{Lemma \stap} There exist constants $C_1,C_2$ such that
$$C_1ILT(x,y)\le d_{PF}(x,y)\le (ILT(x,y)+1)(d_\infty(x,y)+C_2)$$
for all vertices $x,y\in\tg$.
\fthm

\demo{Proof} The first inequality is clear, since an illegal turn
involves at least
one top edge.
For the other, choose a path $\gamma $ from $x$ to $y$ with minimal number of
illegal turns and divide it into legal subpaths $\gamma _j$. We define a subset
$J\inc\gamma $ in the following way: the intersection of $J$ with
$\gamma _j$ is
the maximal subpath
$J_j\inc \gamma _j$ bounded by points of the  full preimage
$r\mi([x,y])$ (where $[x,y]$ denotes the segment  joining $x$ and $y$
in the tree
$\wtilde G$). Note that
$\gamma
\setminus J$ consists of at most $ILT(x,y)
$ intervals.

If $u,v$ bound   a component of   $\gamma \setminus J$, they have
the same image in
$\wtilde G $. Therefore their simplicial distance (hence also their
PF-distance) is bounded. We complete the proof by showing that the
PF-length of a
    $J_j$ is bounded by  the sum of $d_\infty(x,y)$ and a  constant.

Recall that $\pi :\wtilde G \to T$   has
backtracking bounded by some $C'$, and that
$i$ and $\pi \circ r$ are  $C''$-close for some $C''$. If $u,v\in
r\mi([ x,y])$, then in $T$ we have $$d_T(\pi r(u),\pi r(v)) \le d_T(\pi (x),\pi
(y)) +2C'$$ and therefore
$$d_\infty(u,v)\le
d_\infty(x,y)+2C'+4C''.
$$
This bounds the PF-length of a   $J_j$, because
$d_{PF}(u,v)=d_\infty(u,v)$ if
$u,v$ are joined by a legal path.
\cqfd
\enddemo

\head\sect Proof of Theorem \leredu\endhead

We fix $\alpha \in\Aut(F_k)$, and we assume that there is no simplicial $\alpha
$-invariant tree{}  with trivial arc  stabilizers.
Let $G$ be  given  by [\BFHt, 5.1.5], as  above. The
top stratum $G_t$ is
exponentially growing, since otherwise there is a simplicial $\alpha
$-invariant
tree (obtained from $\wtilde G$ by collapsing all zero edges). Let $\G$ and
$T$ be as in Sections 6 and 7.

\subhead \subsect Paired train tracks   \endsubhead


    We define a train track   $G'$ for
$\alpha \mi$ by    applying  [\BFHt, 5.1.5] to $\alpha \mi$ and the
free factor system $\F=\F(G\setminus G_t)$ (consisting of fundamental
groups   of
noncontractible components of   the zero part of $G$).

The top stratum of $G'$ is also exponentially growing (otherwise
there is a  simplicial invariant tree). Since the train track maps
given by [\BFHt, 5.1.5] are
reduced, the free factor systems
$\F$ and
$\F'$ (associated to the zero parts of $G$ and $G'$) are equal (compare the
definition   of pairing in [\BFHt, \S \kern.15em 3.2]).

All the constructions of \S \kern.15em 6 may be applied to $G'$. We use $'$ to
denote the corresponding objects. In particular, we have a graph
$\G'$ (obtained
from $G'$ by adding a shortcut   if needed), a map $f':\tg'\to\tg'$,  a PF
metric on
$\tg'$, and an \Rt{} $T'$.

Since $\F=\F'$, the set of conjugacy classes of $F_k$ represented by
loops in the
zero part is the same for $G$ and $G'$. Furthermore, the fundamental groups of
non-contractible components of the    zero parts of $\G$ and $\G'$
map onto the same
subgroups in
$F_k$ (this follows from Lemma \ellipt{} and the remarks following
it, noting that
$\alpha
$ and
$\alpha
\mi$ have the same invariant conjugacy classes). This implies that
$T$, $T'$ have
the same elliptic elements, and also:

\snom\quasiim
\thm{Lemma \stap} The spaces $\tg\PF$ and $\tg'\PF$ are $F_k$-equivariantly
quasi-isometric.
\fthm

Note that, of course, any equivariant map from $\tg $ to $\tg' $ is a
quasi-isometry for
the {\it simplicial\/} metrics.

\demo{Proof}
    This is standard, and we  only sketch an argument (compare
Proposition 3.1 of [\Far]).
Without changing the quasi-isometry type of $\tg\PF$, or images of
the zero parts
in $F_k$, we may contract to a point  a zero edge of $\G$ with
distinct endpoints, or
contract a top edge if its endpoints are distinct and at most one
touches the zero part. Using these operations and their inverses, we
may assume
that $\G$ has the following
standard form: it has one central vertex
$v$, with top loops $\theta _i$ attached, and top edges
$vv_j$ with a zero component $Z_j$ attached at each $v_j$.

The space $\tg\PF$ is then
quasi-isometric to the Cayley graph of $F_k$ with respect to the infinite
generating system consisting of (the images in $F_k$ of) the $\theta _i$'s and
the whole fundamental groups $\pi _1(Z_j,v_j)$. The space $\tg'\PF$
has a similar
structure, and the two Cayley graphs are quasi-isometric because
one can express each element of one generating system as a word of bounded
length in the other system.
\cqfd\enddemo

\subhead  \subsect  Creating a fixed point
\endsubhead

The final argument  will require  $f$ to have a fixed point   in $\tg$. As
this may not be the case, we have to enlarge $\tg$  (compare [\Luconj,
\S\kern.15em 6]). There are two cases.

If the fixed point $Q$ of $H$ is in $\ov T\setminus T$, there is an eigenray
$\rho $ (see \S \kern.15em 1). Choose
$x\in\wtilde\G$ with $i(x)\in\rho $. By Lemma \iterleg,   there is a legal
path between
$f^p(x)$ and
$f^{p+1}(x)$ for $p$ large. Recall that legal paths map
PF-isometrically into $T$.
We may assume that $y=f^p(x)$ and all its further images by $f$ are
interior points
of  top edges (this rules out countably many choices for the image of $x$ on
$\rho
$).

We now attach an extra edge
$[R,y]$ to $\tg$, and  an edge $[w  R,wy]$ for every $w\in F_k$. We obtain
 a space $\th$ with a cocompact action of $F_k$.  We extend $f$ to $\th$
by mapping $[R,y]$ to the segment
$[R,f(y)]$ (which contains  $y$), keeping the relation $\alpha
(w)f=fw$ satisfied.
The new edges are top edges, with PF-length chosen so that
$d\PF(R,f(y)) =\lambda
d\PF(R, y)$. One of the new turns at $y$ is legal, the other is not.

Everything said above extends to the enlarged space  $\th$. In
particular, Lemmas
\iterleg{} and \legalge{} still hold, and  $\th\PF$ is
quasi-isometric to $\tg\PF$.  The associated
\Rt{} (still denoted by $T$)
consists of the previous $T$ together with the orbit of $Q$. The
image of the new
edge
$[R,y]$ in the tree is the initial segment $[Q,i(y)]$ of the eigenray.

    If
$H$ has a fixed point in $T$, lift it to
$x\in\wtilde\G$ and choose $p$ so that   there is a legal
path between
$y=f^p(x)$ and
$f(y)=f^{p+1}(x)$. Since legal paths map PF-isometrically into $T$,
and
$i\circ f=H\circ i$, this path contains only zero edges.  If $f(y)=y$, we
define $\th=\tg$. If  $f(y)\neq y$, we attach   edges
$[wR,wy]$ as above, but now the new edges are   zero edges.

\subhead \subsect   The main argument
\endsubhead

Fix   train tracks $G,G'$ as in the beginning   of this section. Let
$\G,\G'$ be the corresponding train tracks with shortcuts. Define $\th,\th'$
as just explained, so there are points $R,R'$ in $\th$, $\th'$ fixed by
$f,f'$. They project to the
fixed points
$Q$, $Q'$ of
$H$ and $H'$ in $T,T'$.
We denote by $\C$, $\C'$  the zero components of $\th$, $\th'$
containing $R, R'$.
Recall that $\bo\th$ and $\bo\th'$ are identified to $\bo F_k$.

The stabilizer of $\C$  in $F_k$ is $\st Q$ (see
\S\kern.15em 6). It may be characterized as the only stabilizer which
is invariant
under
$\alpha $ as a subgroup (not just up to conjugacy). In particular,
$\st Q=\st Q'$.

We always
assume that the numbers
$ C,K, C_1,C_2$ introduced in \S\kern.15em 7 work for both $\th$ and $\th'$.
Unless mentioned otherwise, distance refers to the simplicial metric.

We consider $X\in \bo F_k$ with   $Q(X)=Q$.   The numbers $p,\mu ,  \nu $
introduced in the next lemmas will not depend on $X$. These lemmas apply
to the whole $\da$-orbit of $X$ since
$Q(\da^n(X))=Q$ for every
$n\in\Z$.

\snom\contract
\thm{Lemma \stap} Given $\varepsilon  >0$, and $K$ sufficiently
large, there exist
$p=p(\varepsilon )$ and
$\mu =\mu (\varepsilon ,K)$ such that
$$d_{PF}(R,f^p(P))\le
\varepsilon
\ d_{PF}(R,P)+\mu  $$ for all points
$P$ on a K-quasi-geodesic ray $ \gamma=[R,X)\inc\th$  with $Q(X)=Q$.
\fthm

\demo{Proof}
Choose $\delta $ with $\delta (BBT(i)+C_2)<C_1\varepsilon $, where
$C_1,C_2$ come
from Lemma \borneilt. Fix $p$ as in Lemma \iterleg{} (depending on $\delta $).

    Let $\gamma _p$ be a K-PF-geodesic
from $R$ to $\da ^p(X)$. We claim that  $f^p(P)$ is $D$-close to some
$S\in\gamma _p$,  with  $D$ depending only on $K$ and $p$ (but not on $P$
or $X$). This is because $f^p$ is a
quasi-isometry, so $f^p(\gamma )$ is a quasi-geodesic from $R$ to $\da
^p(X)$ with quasi-geodesy constants depending only on $K$ and $p$.

We may assume that $D$ also bounds $ILT(f^p(P),S)$ and $d\PF (f^p(P),S)$. Then
$$ ILT(R,S)-D\le  ILT(R,f^p(P))\le \delta ILT(R,P)\le\delta /C_1\
d_{PF}(R,P)
$$ by    \iterleg{} and \borneilt, and
$$d_{PF} (R,f^p(P))\le d_{PF}(R,S)+D\le
(ILT(R,S)+1)(d_\infty(R,S)+C_2)+D
$$
by \borneilt.

Note that $d_\infty(R,S)$ is the distance between $Q=i(R)$ and $i(S)$ in $T$.
Since
$\gamma _p$ is a K-PF-geodesic with origin $R$, and its point at infinity $\da
^p(X)$ satisfies
$Q(\da ^p(X)) =Q$,  Corollary \coror{} yields $d_\infty(R,S)\le
BBT(i)$.
We obtain an upper bound for $d_{PF} (R,f^p(P))$, which is a linear function of
$d_{PF}(R,P)$ with slope $(\delta /C_1)(BBT(i)+C_2) $.
The  lemma follows.
    \cqfd\enddemo

Lemma \quasiim{} yields an $F_k$-equivariant map $\varphi
:\th\to\th'$ which is a
quasi-isometry for the PF-metrics. It is also a quasi-isometry for
the simplicial
metrics. Similarly, choose $\psi  :\th'\to\th$.   Let $L$ be an upper bound
for all quasi-isometry constants. Fix
$\varepsilon$ with $2L^2\varepsilon <1$ and
choose $p$ as in Lemma
\contract.

\snom\dilat
\thm{Lemma \stap} Given $K'\ge0$, there exists  $\nu \ge0$ such that
$$d_{PF}(R',f'{}^p(P'))\ge2d_{PF}(R',P')-\nu $$   for all points  $P'$ on a
    K'-quasi-geodesic ray $ \gamma'=[R',X)\inc\th'$
with $Q(X)=Q\in  T$.
\fthm

Note that the equality $Q(X)=Q$ does take place  in $T$ (not  in $T'$).

\demo{Proof} Like $\nu $, the numbers $D_1,D_2, D_3$ introduced later
in this proof
will depend on the choices made above, and on $K'$, but not on $P'$
or $X$. There
exists
$D_1$  such that
$\psi (P')$ belongs to a
$D_1$-quasi-geodesic $\theta $ from $R$ to $X$. Choose a $K$-quasi-geodesic
$\gamma
$ from
$R$ to $\alpha ^{-p}(X)$. Since $\theta $ and $f^p(\gamma )$ are
quasi-geodesics
with the same endpoints,
$\psi (P')$ is
$D_2$-close to a point of the form
$f^p(P_1)$, with $P_1\in\gamma $. Finally, we observe that $\varphi (P_1)$ is
$D_3$-close to $f'{}^p(P')$, because $f'{}^p\circ\varphi \circ f^p$ is
$F_k$-equivariant, hence at a bounded distance from $\varphi $.

Now we consider PF-distance. Working modulo additive constants, we have
$$d\PF(R',P')\le Ld\PF(R,f^p(P_1))\le L\varepsilon  d\PF(R,P_1)\le
L^2\varepsilon
d\PF(R',f'{}^p(P'))
$$
(the second inequality comes from Lemma \contract, the other two from
the quasi-isometry
properties of $\varphi $ and $\psi $).  Lemma \dilat{} follows since we have
chosen
$ L^2\varepsilon<\frac12$.
\cqfd\enddemo

Fix $K'$ so that Lemma \geodspec{} applies in $\th'$.

\snom\legaliser
\thm{Lemma \stap} If $Q(X)=Q$, there exists a legal K'-quasi-geodesic
between $R'$
and
$X$ in $\th'$.
\fthm

\demo{Proof} Let $\gamma '$ be a K'-PF-geodesic ray between $R'$ and $ X$. Fix
$P\in\gamma '$. By
\truc{} (and its extension to quasi-geodesic rays), it is $C$-close to all
quasi-geodesics between
$R'$ and
$ X$. We apply this to $f'{}^{np}(\gamma '_{n})$, where $n$ is a
large integer and
    $\gamma '_n$ is a K'-PF-geodesic between $R'$ and $\da ^{np}(X)$.
We see that $P$
is   $C$-close to $f'{}^{np}(P_n)$ for
some $P_n\in\gamma '_{ n}$.

By \dilat{} applied to $\gamma '_n$ we have $$
d_{PF}(R',f'{}^{np}(P_n))-\nu \ge 2^n(d_{PF}(R',P_n)-\nu),$$ so
$d_{PF}(R',P_n) \le \nu +1$  for $n$ large since the left hand side
is bounded. This
gives a uniform bound for
$ILT(R',P_n)$, and by \iterleg{} we deduce
    $ILT(R',f'{}^{np}(P_n))=0$ for $n$ large.

By \geodspec, there exists a
legal K'-quasi-geodesic between $R'$ and
$f'{}^{np}(P_n)$. Since  $f'{}^{np}(P_n)$ is $C$-close to $P$, we conclude by a
diagonal argument, letting $P$ go to infinity on $\gamma '$.
\cqfd\enddemo

Since $Q(\da^n(X))=Q$, we get:

\snom\legaliserr
\thm{Corollary \stap}
If $Q(X)=Q$, then for every $n\in\Z$ there exists a legal
K'-quasi-geodesic between
$R'$ and
$\da^n(X)$ in $\th'$. \cqfd
\fthm

For notational simplicity we state the next results in $\th$, even
though we will
apply them in
$\th'$.

\snom\eig
\example{Remark \stap} If $\rho $ is an eigenray of $H$ (so that
$X=j(\rho )$ is a
fixed point of $\da$), there exists a legal quasi-geodesic $\gamma $
between $R$ and
$X $ in
$\th$ (take
$x\in
\th$ mapping into $\rho $; then $f{}^n(x)\to X\in\bo\th$, and for $n$ large
there is a legal $K$-quasi-geodesic between $R$ and $f{}^n(x)$ by
\iterleg{} and
\geodspec; in fact, $\gamma $ consists of an initial segment contained in the
zero part, followed by  a legal ray in   $\th$). More generally,
there is a legal
quasi-geodesic   between $R$ and $wX $ for $w\in\st Q$.
\endexample

The following lemma is a converse to this remark.

\snom\ptper
\thm {Lemma \stap} Let $X\in\bo F_k$. Suppose that for every $n\in \N$
there exists a
legal
    quasi-geodesic ray between $R$ and $\da^{-n}(X)$ in $\th$. If
$X\notin\bo\st Q$,
there exist $q\ge1$ and
$w\in \st Q$ such that $X$ is a fixed point of  $\bo(
i_w\circ\alpha^q)
$.
\fthm

\demo{Proof} Let $E,E'$ be oriented top edges with origin in $\C$, the zero
component of $\th$ containing $R$.
     Their images in
$T$ are non-degenerate arcs with origin $Q$. In the special situation
that $E,E'$
are in the same
$F_k$-orbit, but distinct,
$i(E)$ and
$i(E')$   don't overlap (because arc stabilizers of $T$ are
trivial). When $E,E'$ are arbitrary (but distinct),  we get a positive
lower bound for possible overlaps between $i(E)$ and
$i(E')$, hence also for possible overlaps of images of
    legal quasi-geodesics with origin in $\C$.

Let $\gamma $ be a legal quasi-geodesic from $R$ to $X$. If
it contains only zero edges, then $X\in\bo\st Q$. From
now on we assume that $\gamma $ has positive PF-length. Its image in
$T$ is a non-degenerate (possibly open) segment
$s(X)$ with origin $Q$, which depends only on $X$ (Lemma \llq{}
implies that $s(X)$
is the set
$B_X$ defined in \S\kern.15em 1, possibly with $Q(X)$ added or
removed).  Note that
$s(\da (X))=H(s(X))$.

We claim that there exist $q\ge1$ and $w\in\st Q$ such that  $s(X)$ and
$s(w\da^q(X))$ have nontrivial overlap.

There is an action of $f_0$ on the   set of oriented top edges of
the finite graph $\h=\th/F_k$ (associate to each top edge the first top edge
of the  image edge path), so there
exists
$r$ such that
$f_0^r$ of every element is periodic.
     Considering the image by $f^r$ of a legal quasi-geodesic from $R$ to
$\da^{-r}(X)$ in $\th$, we see that there exists a legal quasi-geodesic
$\gamma$ from  $R$ to $X$ such that the initial top edges of $ \gamma $ and
some
$f^q(\gamma )$ are in the same $F_k$-orbit (hence in the same $\st
Q$-orbit since
$\st Q$ is the set of deck transformations mapping $\C$  to itself).
This proves
the claim.
 
Let
$\beta =i_w\circ\alpha ^q$, so $w\da^q(X)=\db(X)$. Since $\db
^{-n}(X)$ is in the
same
$\st Q$-orbit as
$\da ^{-nq}(X)$, there is a legal quasi-geodesic $\gamma _{-n}$ from $R$ to
$\db ^{-n}(X)$. If the overlap between $s(X)$ and $s(\db(X))$ is finite,
then   the overlap between $i(\gamma _{-(n+1)})$ and $i(\gamma _{-n})$
is finite and proportional to $\lambda ^{-nq}$.
This is a contradiction since we have seen that it cannot be arbitrarily small.

It follows that $s(X)$ equals $s(\db(X))$ and has infinite length, so
that it is an
eigenray of the homothety $wH^q$ associated to $\beta $. We conclude that
$X=j(s(X))$ is a fixed point of $\db$.
\cqfd\enddemo

We can now  conclude.

\demo{Proof of Theorem \leredu} Recall that $\st Q=\st Q'$.
    If
$Q(X)=Q$, we may  apply    \ptper{} in $\th'$ (thanks to \legaliserr). If
$X\notin\bo\st Q$, we obtain that   $X$ is a fixed point of $\bo(i_w\circ\alpha
^{-q})$ (the exponent is negative because $\G'$ is a train track for
$\alpha \mi$).
Of course this implies that
$X$ is a fixed point of $\bo(i_{\alpha (w^{-q})}\circ\alpha ^q)$, as required.
\cqfd\enddemo

\head\sect More on the dynamics\endhead

\subhead \subsect Products of trees\endsubhead

The techniques used in the previous sections also  give:

\thm{Theorem \stap} Given $\alpha \in\Aut(F_k)$, there exist an
$\alpha $-invariant
\Rt{} $T$ and an $\alpha \mi$-invariant
\Rt{} $T'$,   as in Theorem \arbre, and there exists
$\varepsilon >0$, such that for
every
$g\in F_k$ one of the following holds:
\roster
\item $g$ is elliptic in   $T$ and $T'$;
\item $g$ is hyperbolic in   $T$ and $T'$, and has translation length
$>\varepsilon
$ in $T$ or in $T'$ (or in both).
\endroster
\fthm

This theorem was proved in [\BFHg] and [\Ludis] for $\alpha $ irreducible with
irreducible powers
(no proper
free factor of
$F_k$ is $\alpha $-periodic, up to conjugacy).
It means that the diagonal action of $F_k$ on $T\times T'$ is
discrete. See [\Gui] for further results about such actions.

\demo{Proof} The result is clear if there exists a simplicial $\alpha
$-invariant
tree $T$, with $T'=T$. If not, we let  $\G$, $\G'$, $T$, $T'$ be  as in
\S\kern.15em 8. We use the
same notations.

    We already know that
$T$ and $T'$ have the same elliptic elements. We assume that there is
a sequence
$g_n$ with $\ell(g_n)$ and $\ell'(g_n)$ positive and going to $0$, and we argue
towards a contradiction.

If $x,y\in \tg$ satisfy
$d_\infty(x,y)\le 1$  and $ILT(x,y)\ge 1$, then by \borneilt{}
the ratio between
$d\PF(x,y)$ and
$ILT(x,y )$ lies between $C\mi$ and $C$ for some $C>1$. Similar considerations
    apply to $\G'$,  we fix $C$ working for both $\G$ and $\G'$.

Represent (the conjugacy class of) $g_n$ by a PF-geodesic loop
$\gamma _n\inc\G$.
We write $ILT(\gamma _n)$ for the number of illegal turns of $\gamma _n
$.
    If $ILT(\gamma _n)$ remains bounded, then by \iterleg{} there exists
$r$ such that $\alpha ^r(g_n)$ is represented by a legal loop. The PF-length of
that loop is bounded away from $0$, and  so is
$\ell(g_n)=\lambda ^{-r}\ell(\alpha ^r(g_n))$  by \legalge.
Assume therefore that $ILT(\gamma _n)$  goes to
infinity.

Let $L$ be as above (quasi-isometry constant). Fix $\delta >0$ with
$LC^2\delta <1$.
Let $p$ be given by \iterleg{} (applied to $\delta $ in both $\tg$ and $\tg'$).
Choose $x\in\tg$ projecting into  $\gamma _n$ in $\G$ and into
     the axis of $g_n$ in $T$, and let $y=g_nx$. Note that
$d_\infty(x,y)=\ell(g _n)$
and
$d_\infty(f^p(x),f^p(y))=\lambda ^p d_\infty(x,y)$
go  to
$0$ as $n\to \infty$.

For $n$ large we have
$$d\PF(f^p(x),f^p(y))\le C \cdot ILT(f^p(x),f^p(y))\le C\delta  \cdot
ILT(x,y)\le
C^2\delta  d\PF(x,y),
$$
showing that
$h_n=\alpha ^p(g_n)$   may  be represented   in
$\G$ by a  loop of PF-length less than $ C^2\delta |\gamma _n|\PF$,
hence in $\G'$
by a loop of PF-length less than $LC^2\delta |\gamma _n|\PF+L$.

Since $\ell'(h_n)$
goes to $0$, we may apply the same argument to $h_n$ in $\tg'$, and we get
    $|\gamma _n|\PF\le (LC^2\delta )^2|\gamma _n|\PF+D$ for some constant
$D=D(L,C,\delta )$.  This shows that
$|\gamma _n|\PF$, hence $ILT(\gamma _n)$,  is bounded, a
contradiction.
\cqfd\enddemo

\subhead \subsect Dynamics of irreducible automorphisms\endsubhead

Consider $\alpha \in\Aut(F_k)$. For simplicity   assume   that
all periodic points of $\ov\alpha $ are fixed points (this may be
achieved by raising
$\alpha $  to some power). Recall [\GJLL] that fixed  points of $ \da$ not in
$\bo\fix\alpha
$ are either attracting or repelling, and the action of $\fix\alpha $ on
$\fix\da\setminus\bo\fix\alpha
$ has finitely many
    orbits.

When $\fix\alpha $ is trivial, $\fix\da$ is   the vertex set
of a  {\it finite bipartite graph $\Gamma $\/}, with an edge from a
repelling point
$X_1$ to an attracting point $X_2$  if and only if there exists
$X\in\bo F_k$ with
$ \lim_{n\to+\infty}\da ^{- n}(X)=X_1$ and $
    \lim_{n\to+\infty}\da ^{  n}(X)=X_2$. Note that
every component of $\Gamma $
contains at least two vertices.

Recall that an automorphism $\alpha $ is irreducible with irreducible
powers (iwip) if no
proper free factor of
$F_k$ is $\alpha $-periodic (up to conjugacy).

     Dynamics of geometric iwip automorphisms (induced by a
pseudo-Anosov homeomorphism
of a compact surface with one boundary component) is well-understood.
Because there
is an invariant cyclic ordering on $\bo F_k$, the graph $\Gamma $ (defined when
$\fix\alpha $ is trivial) is either a single edge or it is
homeomorphic to a circle (the
second author has conjectured that this property leads to a
characterization of geometric
automorphisms).

We focus on   non-geometric
iwip automorphisms (in this case, $\fix\alpha $ is always trivial [\BH]).
We do not know which graphs $\Gamma $ may appear in this context. We
only prove:

\thm{Theorem \stap} Let $\alpha \in\Aut (F_k)$ be irreducible with irreducible
powers, not geometric. Assume that all periodic points of $\da$ are \fp s. Then
either the graph
$\Gamma
$ has exactly two vertices, or every component of $\Gamma $ contains
strictly more
than two vertices.
\fthm

\demo{Proof}
Let $T$ be an invariant \Rt{}     as in Theorem \arbre{} (it is unique up to
rescaling in this case, see [\LLquat]). It is well-known that
$\fix\alpha $ is trivial, and the action of $F_k$ on $T$ is free (this may be
deduced from [\BH] and Lemma \ellipt). We have $\lambda >1$, and we let
$Q$ be the
\fp{} of $H$ (in $\ov T$).

By Proposition \bij{}, all components of $\ov
T\setminus\{Q\}$ are fixed by $H$.
If $g\in F_k$ is nontrivial, then   $\alpha ^n(g)$ converges to the
attracting \fp {} $j(\rho )$, where
$\rho $ is the eigenray of $H$ contained in the
    component   containing $gQ$ (see the proof of Theorem \fevr). There is a
similar relation between the
$\alpha
\mi$-invariant tree $T'$ and repelling \fp s of $\da$.

We may assume that $\da$ has at least two attracting \fp s,
and two repelling ones (otherwise, the result is  trivial). By   \bij{}
this implies that $\ov
T\setminus\{Q\}$ has at least two components, so  $Q$ belongs to
$T$ (not just to
$\ov T$). Similarly $Q'\in T'$.
We now argue by way of contradiction, assuming that some
component of $\Gamma $ has only two vertices, an attracting fixed
point $X $ and a
repelling point $ X'$. Let
$\C,
\C'$ be the corresponding components of $
T\setminus\{Q\}$, $
T'\setminus\{Q'\}$ respectively, and $\rho $ the eigenray of $H$ contained in
$\C$.
 
Viewing nontrivial elements of $F_k$ as hyperbolic isometries of
$T$, we claim that {\it there exists $g\in F_k$ whose translation axis $A_g$
passes through $Q$ and intersects $\rho $ in a segment strictly longer than the
translation length
$\ell(g)$ (i.e\.
$gQ$ is an interior point of
$A_g\cap
\rho
$). }

Assuming this claim temporarily, we complete the proof as follows.
The point $g\mi Q$ belongs to a component $\C_1$ of $
T\setminus\{Q\}$ distinct from $\C$. Since the edge $X'X$ of $\Gamma $ is
isolated, the point $g\mi Q'\in T'$ belongs to a component $\C'_1$ of $
T'\setminus\{Q'\}$ distinct from $\C'$ (otherwise $\Gamma $ would contain an
edge between $X'$ and the attracting \fp {} corresponding to $\C_1$).
Now consider
$g\mi\alpha ^p(g)$, for $p$ large.  The point $g\mi\alpha ^p(g)Q=g\mi
H^p(gQ)\in T$ belongs to $\C$ because of our choice of $g$.
    In $T'$, on the other  hand, the point $g\mi\alpha
^p(g)Q'=g\mi (H')^{-p}(gQ')$ is close to $g\mi Q'$ and therefore   belongs to
$\C'_1$. It follows that $X $ is joined to  two distinct vertices of
$\Gamma $, a
contradiction.

There remains to prove the claim. We use the terminology of  \S\kern.15em 6.
Since $\da$ has at least two attracting \fp s, it follows from
[\GJLL, p\.~431] that
$ f$ has a \fp{} $R\in\tg$. By Remark \eig, there exists a legal quasi-geodesic
$\gamma
$ between $R$ and $X=j(\rho )$.  By irreducibility of $\alpha $, the
projection of $\gamma  $ onto $\G $ passes again over its initial top
edge (with the
same orientation). This defines a loop  in $\G $, and a
nontrivial
$g\in F_k$ mapping an initial segment of $\gamma $ into $\gamma $ (in an
orientation-preserving way). This is the required $g$. \cqfd
\enddemo

\subhead \subsect The number of periods \endsubhead

Given $\alpha \in\Aut(F_k)$, recall [\LLdeux] that
    periods of elements of $  F_k$ for $ \alpha $ are bounded
by $A_k$, the maximum order of torsion elements in $\Aut(F_k)$, and periods of
elements of $\bo F_k$ are bounded by $M_k=2kA_k$. For $k$ large, one has
$\log A_k\sim\log M_k\sim\sqrt{k\log k}$.

Examples of automorphisms with  many periods may be
     constructed as follows. Let $p$ be a prime number. Let $\sigma $ be
a permutation
consisting of one cycle of order $p'$ for each prime number $p'\le
p$. It defines
a periodic automorphism $\alpha $ of $F_k$, with $k=2+3+\dots+p$. The
periods of
$\ov\alpha
$ are exactly the divisors of $ 2\cdot3\cdots p$. A simple computation (based
on the prime number theorem) shows that, for
$p$ large, the logarithm of the number of periods of $\ov\alpha $ is
asymptotic to
$2\log2\cdot
\sqrt{\frac k {\log k}}$.

\thm{Theorem \stap}
Given $\alpha \in\Aut(F_k)$,   the number  of periods of periodic
points of $\ov\alpha $ is bounded by $N_k$, with $\log
N_k\sim2\log2\cdot \sqrt{\frac k {\log k}}$.
\fthm

\demo{Proof}  First consider periods of elements of $F_k$. Let
$P(\alpha )$ be the
periodic subgroup (consisting of all $\alpha $-periodic $g\in F_k)$.
Since $P(\alpha
)$ is the fixed subgroup of some power of $\alpha $, it has rank at most $k$ by
[\BH]. The number of
periods of $\alpha $ is   bounded by the number of divisors of
$o(\alpha )$, where $o(\alpha )$ is
    the order  of $\alpha $ in $\Aut(P(\alpha ))$. Now $o(\alpha )\le A_k$, and
    the number $d(n) $ of divisors of $n$ satisfies
$\ds\limsup_{n\to+\infty}\frac{\log d(n) \log\log n}{\log n}=\log2$ (see [\HW,
Theorem 317]). Replacing
$n$ by $A_k$ with $\log A_k \sim\sqrt{k\log k}$ gives
$N_k$ as in the theorem.
 
If $X \in\bo F_k$ is periodic, and no $g\in F_k$ has the same period,
the proof of
Theorem 2.1 of [\LLdeux] shows that the period of $X$ divides $rs$, where $r\le
2k$ and $s$ is the period of some $g\in F_k$. The estimate therefore
also holds for
the periods of $\da$ since $\log 2kN_k\sim\log N_k$.
\cqfd\enddemo

\subhead \subsect Automorphisms with many fixed points \endsubhead

We give a short proof of a result  due to Bestvina, Feighn,
Handel [\BFHs], improving their lower bound  from 3 to 4.

\snom\bfh
\thm{Proposition \stap} For any outer automorphism $\Phi $ of $F_k$, $k\ge2$,
there exist
$q\ge1$ and $\beta  \in \Aut(F_k)$ representing $\Phi ^q$ such that
$\db$ has at least
four
\fp s. If $u\in F_k$ is fixed by some
$\alpha \in\Aut(F_k)$ representing $\Phi $, we may require $\beta (u)=u$.
\fthm

\remark{Remark}
    With the terminology  of [\GJLL, \S\kern.15em 6], we shall prove that any
$\Phi \in\Out
(F_k)$ has a power with positive index. It is not always possible to
take $q=1$: if
$\alpha \in\Phi $ cyclically permutes the elements of a free basis of $F_k$, it
follows from [\CL] that $q$ has to be divisible by~$k$.
\endremark

\demo{Proof} Starting with   $\alpha \in\Phi $, we will keep
replacing it by a power
$\alpha ^s$, or by
$i_w\circ\alpha
$,  so as to finally obtain an automorphism (still denoted by $\alpha
$) with  at least
four
\fp s on $\bo F_k$.

Let $T$ be an
$\alpha
$-invariant
\Rt{} as in Theorem \arbre. Replacing $\alpha $ by
$\alpha ^s$ or by
$i_w\circ\alpha
$  amounts to replacing $H$ by   $H^s$
or   $wH$.

By assertion (d) of Theorem \arbre{}, we may assume that $H$ fixes a branch
point
$Q\in T$ and acts trivially on the set of $\st Q$@-orbits in $\pi
_0(T\setminus\{Q\})$. We distinguish several cases.

If $\lambda >1$ and $\st Q$
is trivial, then $\da$ has at least three attracting
\fp s (see Proposition \bij). There is a fourth, repelling, point.
If $\lambda >1$ and $\st Q=\Z$, some $wH$ with $w\in\st Q$ has an eigenray and
there are infinitely many attracting periodic points (because $\fix\da^2$ is
invariant
     under the action of $\st Q$). If $\lambda >1$ and $\st Q$ has rank
$\ge2$, we use
induction on
$k$.

    If
$\lambda =1$, we may assume that $H$ fixes an edge $e=[Q,R]$, and we
collapse every
edge not in the orbit of $e$. In the new tree, $\st Q$ and $\st R$
are non-trivial
and $\alpha $-invariant. This proves the first part of the proposition.

Now assume $\alpha $ fixes some nontrivial $u\in F_k$ (and therefore
$\da$ has two
\fp s $u^{\pm\infty}$). Note that $H$ commutes with $u$. If $u$ fixes
a (unique)
point
$Q\in T$ (in particular if
$\lambda >1$), this point is also fixed by $H$ and we argue
as before.
Finally, suppose $\lambda =1$ and $u$ is hyperbolic. In
this case we may assume that $H$ equals the identity on the axis $A$ of $u$
(replace $H$ by some $u^rH^s$). After possibly
collapsing we get
$Q\in A$ with nontrivial stabilizer and we consider $\alpha _{\mid\st Q}$.
\cqfd\enddemo

\head \sect  Hyperbolic groups \endhead

\subhead  \subsect  Bounding   periods \endsubhead

Let $\Gamma $ be a torsion-free hyperbolic group. Given $\alpha
\in\Aut(\Gamma )$,
let $\fix\alpha $ be its fixed subgroup, and
$P(\alpha )=\cup_n\fix\alpha ^n$ be its periodic subgroup. A subgroup
of $\Gamma $
is a {\it fixed subgroup\/} (resp\. a {\it periodic subgroup\/}) if it equals
$\fix\alpha $ (resp\.
$P(\alpha )$) for some
$\alpha \in\Aut(\Gamma )$.

\thm{Proposition \stap} Every periodic subgroup is hyperbolic. Up to
isomorphism,
there are only finitely many periodic subgroups in a given $\Gamma $.
\fthm

\demo{Proof}  Most arguments come from [\Sh]. There are two cases.

$\bullet$ Suppose $\Gamma $ is one-ended (= freely indecomposable).
By [\Se, Theorems 3.2 and 4.1], the group $P(\alpha
)$, if not trivial or cyclic, is a vertex group in some splitting of
$\Gamma $ with
cyclic edge groups. By [\Guia],   such a splitting is obtained from
the JSJ splitting constructed in [\Bow] by blowing up vertices  with
group equal to
$\Z$, blowing up  surface vertices along disjoint simple closed
curves, and then
collapsing edges. Finiteness of the set of curves on a compact surface (up to
homeomorphism) implies that there are only finitely many possibilities for
$P(\alpha )$ (up to an automorphism of
$\Gamma
$). Furthermore, $P(\alpha )$ is hyperbolic because it is a vertex group in a
splitting with quasiconvex edge groups [\Bow].

$\bullet$ Now suppose that $\Gamma $ is the free product of cyclic groups and
one-ended groups $\Gamma _i$. By the Kurosh subgroup theorem,  $P(\alpha )$
    is
the free product of cyclic groups  and  subgroups
$H_j$ of  conjugates $K_j$ of the
$\Gamma _i$'s. The number of factors is the {\it Kurosh rank\/} of
$P(\alpha )$. It   is finite  because $P(\alpha )$ is the increasing
union of the $\fix(\alpha
^{n!})$, whose Kurosh rank is uniformly bounded    [\CT].

If $\alpha ^p(\Gamma _i)$ meets
$\Gamma _i$ non-trivially, then $\alpha ^p(\Gamma _i)=\Gamma _i$. It
follows that
each
    $K_j$ is $\alpha $-periodic, and therefore $H_j=P(\alpha )\cap K_j$ is the
periodic subgroup of some automorphism of $K_j$. By the first case, $H_j$ is
hyperbolic and has only finitely many possibilities (up to
isomorphism). The same is
therefore true of
$P(\alpha )$.
\cqfd\enddemo

Since $P(\alpha )$ is finitely generated, the restriction of $\alpha
$ to $P(\alpha
)$ has finite order. In particular, every periodic subgroup is a
fixed subgroup and
is quasiconvex [\Ne].

Recall  [\Lev] that for
$P$ torsion-free hyperbolic there are only finitely many conjugacy classes of
torsion elements in
$\Aut(P)$. Conjugate automorphisms having isomorphic fixed subgroups,
we deduce:

\thm{Corollary \stap{} (Shor [\Sh])} Up to isomorphism, there
are only finitely many fixed subgroups in a given torsion-free hyperbolic group
$\Gamma $. \cqfd
\fthm

    We  have also proved part of assertion (1) of Theorem \pth{}:

\snom\bor
\thm{Corollary \stap} Given a torsion-free hyperbolic group $\Gamma $,
there exists $M$ such that, if $g\in \Gamma $ is periodic under $\alpha \in\Aut
(\Gamma )$, then the period of $g$ is at most $M$.
\fthm

\demo{Proof} The period of $g$ divides the order of $\alpha $ in
$\Aut(P(\alpha ))$.
\cqfd\enddemo

\subhead  \subsect  One-ended groups\endsubhead

Theorem \pth{} will first be proved for torsion-free groups, in this subsection
(one-ended groups) and in the next (groups with infinitely
many ends). Groups with torsion will then be considered.

Let
$\Gamma
$ be a torsion-free one-ended hyperbolic group, and
$\alpha
\in\Aut(\Gamma )$. Using Corollary \bor, we shall prove that {\it
every $X\in\bo
\Gamma
$ is asymptotically periodic\/}, with a uniform bound (depending only
on $\Gamma $)
for the period of the limiting orbit.  This will prove
assertions (1) and (3) of Theorem \pth{} for $\Gamma $. Assertion (2) is proved
by similar arguments, or deduced from   (3) since the sequences
$\alpha ^n(g)$ and
$\da^n(g^\infty)$ have the same
limit points when $g$ is not $\alpha $-periodic (see [\LLun, proof of 2.3]).

    We use the
JSJ-splitting of
$\Gamma $, first introduced by Sela. We prefer to follow Bowditch's approach
[\Bow] because of its strong uniqueness properties.
The JSJ-splitting decomposes $\Gamma $ as the fundamental group of a
finite graph of
groups, with  associated
Bass-Serre tree $T$. Since $T$ is constructed   purely from the
topology of $\dg$, the group $\Aut (\Gamma )$ acts on $T$ in the same way as on
$\dg$. In particular, there is an isometry   $H :T\to T$ as in Theorem \arbre.

Edge stabilizers are cyclic. Vertex stabilizers $\st Q$ are
quasiconvex [\Bow],
hence  hyperbolic.
    Furthermore the boundary of $\Gamma $ is the
disjoint union of the set of ends of $T$ (embedded into $\bo \Gamma $
by a map $j$
as in \S\kern.15em 1.d) and the (non-disjoint) union of the boundaries
$\bo \st Q$
of the vertex stabilizers (see [\Bow]).

A vertex stabilizer  $\st Q$ is
    cyclic, or   free (``hanging fuchsian''), or   ``relatively
rigid'' (the subgroup of
$\Out (\st Q)$ consisting of outer automorphisms fixing
stabilizers of edges incident to $Q$ (up to conjugacy) is  finite; as
explained in
[\Lev], this follows from [\BF] and [\Pau]).

If an edge $e$ of $T$ is fixed by $H$, its
stabilizer is an $\alpha $@-invariant cyclic subgroup of $\Gamma $.
By Corollary
\bor,  we may raise $\alpha $ to a
fixed power (depending only on
$\Gamma $) to ensure that any $H$@-periodic edge is in fact   fixed.

Suppose that a vertex $Q\in T$ is fixed by $H$. Then $\st Q$ is $\alpha
$@-invariant, and the induced automorphism has finite order in
$\Out(\st Q)$ (if
$\st Q$ is cyclic or relatively rigid) or is a geometric automorphism of a free
group (if $\st Q$ is hanging fuchsian). It follows that  any
$X\in\bo(\st Q)$
is  asymptotically periodic,
with a uniform bound on the
period.

The proof for arbitrary $X  $   now  is  fairly similar to the proof of Theorem
\simpli. Lemma \maine{} extends to     actions of  torsion-free hyperbolic
groups with trivial or cyclic edge stabilizers (compare [\Dah, \Ka]).
Lemma \main{}
also extends to torsion-free hyperbolic groups.

Let $X\in\bo \Gamma $.
If the isometry $H:T\to T$ is hyperbolic, with axis $A$, then either
$X=j(A^-)$ is
fixed by
$\da$, or $\da^n(X)$ converges to $j(A^+)$ as $n\to+\infty$. If $H$
is elliptic,
let $\P$ be its fixed subtree (equal to the set of periodic points).

First suppose $X\in\bo\st R$ for some vertex $R$. Using remarks made
above, we may
assume that
$R$ cannot be chosen in $\P$. Let $Q$ be the point of $\P$ closest to
$R$. Let $\C$
be the component of $T\setminus\{Q\}$ containing $R$. If $g_p\in\st
R$ converges to
$X$, we have $g_pQ\in \C$ for $p$ large, because otherwise $g_p\in\st Q$ and we
could have chosen
$R$ in $\P$. Now apply Lemmas \maine{} and \main{} as in the proof of Theorem
\simpli. Note that  Remark \sousgr{} provides a uniform bound for the
cardinality of
$\omega (X)$.

If $X=j(\rho )$ for some end $\rho $ of $T$, then either $\rho $ is
an end of $\P$
(and
$X$ is fixed by $\da$), or we apply Lemmas \maine{} and \main{} using   the
point
$Q\in\P$ closest to $\rho $ and the component  $\C$  of
$T\setminus\{Q\}$ containing
$\rho $.

\subhead  \subsect  Free products\endsubhead

Let  $\Gamma $ be torsion-free, with infinitely many ends. We write
$\Gamma =F_k*\Gamma _1*\dots*\Gamma _m$, with each $\Gamma _i$
one-ended. We will
use the invariant
\Rt{}   given by the following result.

\snom\GJLLfreepr
\thm{Theorem \stap} Given $\alpha \in\Aut (\Gamma )$, there exist an
\Rt{} $T$ and a
homothety $H$ satisfying  conditions (a), (b), (c) of Theorem
\arbre{}. Furthermore
each
$\Gamma _i$ ($1\le i\le m$) fixes a point of $T$, and there exists a $\Gamma
$-equivariant injection $j:\bo T\to\bo \Gamma $ satisfying $\da\circ
j=j\circ H$.
\fthm

\demo{Proof}
The construction of $T$ in the case of $F_k$ has been sketched in
\S\kern.15em 6
(see [\GJLL] for details):  equip
$\wtilde G$ with the PF-metric $d\PF$, and consider the metric space
associated to
$d_\infty$ (when $\lambda =1$, simply collapse   components of the zero set of
$\wtilde G$ to points). The proof in the general case is similar, using the
``efficient representatives'' constructed by Collins-Turner in
    [\CT]. The only difference is that the zero set is a 2-complex (not
necessarily a
graph),  but this
difference is irrelevant  as each component gets   collapsed to a point.

Conditions (a) and
(b) of Theorem \arbre{} are proved as in
[\GJLL] (Lemmas 2.7 and 2.8), and $\Gamma _i$ fixes a point in $T$
because it fixes
a component of the zero set in $\wtilde G$. The map $j$ is constructed as in
    [\GJLL, Lemmas 3.4 and 3.5], but we need to show that equivariant
maps from Cayley
graphs of $\Gamma $ to $T$ have bounded backtracking (in the sense of
\S\kern.15em
1.d).

This is true for maps to  $\wtilde G$ (equipped with a simplicial
metric), because
$\wtilde G$ is quasi-isometric to $\Gamma $. Therefore it is also
true for maps to
the tree
$T_0$ obtained from $(\wtilde G, d\PF)$ by collapsing the zero set.
When $\lambda >1$, we further observe that
the map $\ov f:T_0 \to T_0$ induced by $f$ is $\lambda
$@-Lipschitz for the PF-metric, and has backtracking bounded by some
constant $K$.
Letting
$d_p(x,y)= \lambda ^{-p}d\PF(\ov f^p(x),\ov f^p(y))$,
     we deduce
that the identity map from $(T_0,d\PF)$ to $(T_0,d_p)$ has
backtracking bounded by
$K(\lambda
^{-1}+\dots+\lambda ^{-p})$. The canonical map from $(T_0,d\PF)$ to
$T=(T_0,d_\infty)$ is
$1$-Lipschitz, and has bounded backtracking because the series $\lambda
^{-1}+\dots+\lambda ^{-p}+\dots$ converges. It
follows that equivariant maps from Cayley graphs of $\Gamma $ to $T$
have bounded
backtracking.\hfill\Square
 \enddemo

    Recall that the {\it rank\/}
$\rk(J)$ of a group $J$ is the minimum cardinality of a generating
set (not to be
confused with the Kurosh rank  used in [\CT]).

\snom\damiensformula
\thm{Theorem \stap{} (Gaboriau [\Gab])} Let $\Gamma $ and $T$ be as above.
For any $Q\in T$, the stabilizer $\st Q$ has rank $\rk(\st Q)\le
\rk(\Gamma )-1$,
and the action
of $\st Q$ on $\pi _0(T\setminus\{Q\})$ has at most $2\rk(\Gamma ) $
orbits. \cqfd
\fthm

We can now prove Theorem \pth{} for $\Gamma $.
Let $T $ be as in Theorem \GJLLfreepr.
Note that for $Q\in T$ the intersection of $\st Q$ with a conjugate
$g\Gamma _ig\mi$
is either trivial or the whole of $g\Gamma _ig\mi$ (because arc stabilizers are
trivial and $\Gamma _i$ fixes  a point). Thus vertex stabilizers are
free products,
with each factor free or isomorphic to some
$\Gamma _i$. They are
    ``simpler'' than $\Gamma $ by Theorem
\damiensformula{},   quasiconvex (see e.g\. [\Ka, \Swa]), and up to
isomorphism they
belong to some finite set depending only on $\Gamma $.

The proof of assertion (2) of Theorem \pth{} now  proceeds exactly as in the
case of
$F_k$ (Theorem
\fevr), by induction on rank, since the result is already known for one-ended
groups.

To prove assertion (1), we need to bound the period of a
$\da$@-periodic $X\in\bo
\Gamma $.  If $X\in\bo P(\alpha )$, the period is bounded by the
order of $\alpha $
in $\Aut(P(\alpha ))$. If not, then $X$ is attracting or repelling,
hence  belongs
to the
$\omega$-limit set (for
$\alpha $ or $\alpha \mi$) of every
    $g\in
\Gamma
$ close enough to $X$ in $\ov \Gamma $ (see the discussion in \S\kern.15em 4 of
[\LLdeux]). We therefore reduce to controlling the cardinality of $\omega (g)$.

    The arguments from the proof  of Theorem \fevr{} do not provide
uniform bounds, for two reasons.  If $\lambda >1$ and the component
$\C$ is $H$@-periodic, we do not have a bound for the period. If
$\lambda =1$ and
the isometry $H$ is elliptic, we do not have a bound for the period
of its periodic
points.

We first show that,  {\it if $Q$ is a fixed point of $H$, there is a
bound depending
only on
$\Gamma
$ for the period $p$ of an $H$@-periodic component  $\C$ of
$T\setminus\{Q\}$.} By
Theorem
\damiensformula{} we may assume that $H$ acts trivially on the set of
orbits of the
action of $\st Q$ on
$\pi _0(T\setminus\{Q\})$. We then have $H\C=w\C$ for some $w\in\st Q$, and
$H^p\C=w_p\C$ with $w_p=\alpha ^{p-1}(w)\dots\alpha (w)w$. From
$w_p=1$ we get
$\alpha ^p(w)=w$, hence $\alpha ^r(w)=w$ for some $r$ depending only
on $\Gamma $
by Corollary \bor, and finally
$w_r=1$ because $(w_r)^p=w_{pr}=1$. This implies $p\le r$.

Recall that we want to bound the period of $X\in\omega (g)$. If $\lambda >1$,
the arguments from the proof of  Theorem \fevr{} show that  either
$X=j(\rho )$ for some
$H$@-periodic ray, and we are done, or $\omega (g)\inc\bo\st Q$ and we can use
induction on the rank of $\Gamma $.

Now suppose $\lambda =1$.
There is a problem only if
$H$ is elliptic and has periodic points with large periods. It then has a fixed
point $Q$, and using the fact proved above about periodic components  of
$T\setminus\{Q\}$ we may assume that
$H$ fixes an edge
$e$ (this involves raising $\alpha $ and $H$ to some  power depending only on
$\Gamma $). Consider the tree
$T'$ obtained by collapsing all edges of $T$ not in the $\Gamma
$@-orbit of $e$.

It
satisfies the conditions of Theorem \GJLLfreepr{}, and the quotient
graph
$ T'/\Gamma $ has exactly one edge (it is a segment or a loop). We can now
complete the proof by induction, using the special form
of
$\alpha
$  as in [\GJLL, pp\. 442-443] (if for instance
$T'/\Gamma $ is a segment, then
$\alpha $ preserves a nontrivial decomposition of $\Gamma $ as a free
product, and
any attracting periodic point $X$ may be written $X=gX'$ where $g\in
\Gamma $ is
$\alpha $-periodic and
$X'$ is contained in the boundary of a free factor). This completes
the proof of
assertion (1) of Theorem \pth.

The proof of assertion (3)   is the same as for free groups (proof
of Theorem \simpli). Lemma
\finiper{} extends, replacing $M_k$ by the bound obtained above. This
completes the
proof of Theorem \pth{} for torsion-free groups.

\subhead \subsect   Groups with torsion \endsubhead

The goal of this subsection is to extend Theorem \pth{} to virtually
torsion-free
groups. We start with general considerations.

Suppose $\Delta $ has finite index in
a group
$\Gamma $. Then there exists
$k$ such that, for every $\alpha
\in\Aut (\Gamma )$, every coset of $\Gamma $ modulo $\Delta $ is
mapped to itself by
$\alpha ^k$. In particular, $\alpha ^k(\Delta )=\Delta $ for every
$\alpha
\in\Aut (\Gamma )$. Of course $\alpha ^k{}_{|\Delta }$ is
polynomially growing if
$\alpha
$ is.

Suppose furthermore that there is a bound for the order of torsion in
$\Delta $.
Then    {\it there is a bound for periods of elements of $\Gamma $ under
automorphisms of
$\Gamma
$ if there is one for periods of elements of $\Delta $ under automorphisms of
$\Delta $.} To prove this, suppose $\alpha ^p(g)=g$ with $g\in \Gamma
$ and $\alpha
\in\Aut (\Gamma )$. Replacing
$\alpha
$ by $\alpha ^k$, we may assume
$\alpha (g)=wg$ with $w\in \Delta $. Then $\alpha ^p(g)=w_pg$ with
$w_p=\alpha ^{p-1}(w)\dots\alpha (w)w$. As in the previous section,   $w_p=1$
implies
$\alpha ^p(w)=w$, and therefore $\alpha ^r(w)=w$ for some $r$  which
can be bounded
in terms of
$\Gamma $ and
$\Delta $ only. We then write $(w_r)^p=w_{pr}=1$, and we bound the
period of $g$ by
$r$ times the order of $w_r$.

Now suppose that $\Gamma $ is hyperbolic and $\Delta $ is a
torsion-free subgroup
of  finite index. We have just proved assertion (1) of Theorem \pth{} for
periodic orbits of
$\alpha $.
    Since $\Delta $
and $\Gamma $ have the same boundary, assertions (1) (for orbits of
$\da$) and (3)
hold.

To prove assertion (2), consider
$g\in \Gamma $ and $\alpha \in\Aut (\Gamma )$. Replacing $\alpha $ by
$\alpha ^k$,
we may assume $\alpha (\Delta )=\Delta $ and $\alpha (g)=hg$ with
$h\in \Delta $.
Then
$\alpha ^n(g)=h_ng$ with
$h_n=\alpha ^{n-1}(h)\dots\alpha (h)h$, and we conclude by Corollary \trick{}
(proved in $\Delta $ just like in $F_k$).

We also show:

\thm{Proposition \stap}  If $\alpha \in\Aut(\Gamma )$, with $\Gamma $
an infinite,
virtually  torsion-free,
hyperbolic group, then $\da$ has at least two periodic
points. If $\da$ has only one periodic orbit, then this orbit  is the
boundary of
an $\alpha $-invariant virtually cyclic subgroup.
\fthm

\demo{Proof} As in [\LLdeux, proof of 1.1]. If $P(\alpha )$ is
finite, assertion
(2) of Theorem \pth{} provides both an attracting periodic orbit and a
repelling one. If
$P(\alpha )$ is virtually
$\Z$, its boundary gives two fixed points, or an orbit of order 2. If
$P(\alpha )$
is non-elementary, there are uncountably many periodic orbits.
\cqfd\enddemo

\head  {\sect Examples  and questions}\endhead

\subhead \subsect Examples \endsubhead

Let $\alpha \in\Aut(F_k)$.  
Fixed points of $\da$ not in $\bo \fix \alpha $ are either attracting
or repelling [\GJLL].
Now consider the automorphism $\alpha $
of $F_2$ mapping $a$ to $a$ and $b$ to $aba$.  The group
$\fix\alpha =\langle a\rangle$ is infinite cyclic, and its two limit points
$a^{\pm\infty}$ are isolated fixed points of $\da$ which are
{\it half-attracting and
half-repelling:\/} if $X\in\bo F_2$ is not $a^{\pm\infty}$, then the
limit of $\da^n(X)$ as
$n\to+\infty$ is either $a^\infty$ or $a^{-\infty}$, depending on whether
the first occurrence
of $b^{\pm1}$ in $X$ is $b$ or $b\mi$. 
The automorphism $\alpha $ is the square of $\beta  :a\mapsto a\mi, 
b\mapsto a\mi
b\mi$. The set $\{ a^\infty,a^{-\infty}\}$ is a half-attracting,
half-repelling orbit of  period two of $\db$. 

The following example of a {\it parabolic orbit\/} is due  to A. Hilion
[\Hi]. Define $\gamma $ on $F_4$ by $a\mapsto a$, $b\mapsto ba$, $c\mapsto
ca^2$,
$d\mapsto dca$. For $g=bad\mi$, both sequences $\gamma ^n(g)$ and $\gamma
^{-n}(g)$ converge to $ba^{-\infty}$ as $n\to+\infty$.
 
Here is an example where 
    $\fix\alpha $ is cyclic, and its limit
points are not isolated
as fixed points of $\da$. Consider a homeomorphism $h$ of a compact
surface $\Sigma $
fixing a separating simple closed curve $C$ pointwise. Assume that
$h$ induces a
pseudo-Anosov homeomorphism on each of the complementary subsurfaces $\Sigma
_\ell$ and $\Sigma _r$, and $h$ twists non-trivially around $C$ in  $\Sigma
_\ell$ (but not in $\Sigma _r$). Let $\alpha $ be the automorphism
induced on $\pi
_1(\Sigma )$ (with basepoint on $C$), and $\da$ the homeomorphism
induced on its
boundary. Note that the boundary is a circle (if $\Sigma $ is closed)
or a cyclically
ordered Cantor set (if $\Sigma $ has a boundary). The map $\da $ has
two fixed points
associated to the invariant cyclic subgroup $\pi _1(C)$. They divide
the boundary into
two intervals $I_\ell$ and $I_r$. On $I_\ell$, the map $\da$ has no
fixed  point (because of
the non-trivial twist). On
$I_r$, there are infinitely many attracting fixed points, and
infinitely many repelling ones (but only finitely many orbits under
the action of $\pi
_1(C)$; they correspond to singular leaves of the invariant
foliations of $h_{|\Sigma
_r}$ issuing from singularities belonging to $C$). They alternate on
$I_r$, and accumulate onto both endpoints of
$I_r$.

\subhead \subsect Free groups\endsubhead

Consider $\alpha \in\Aut(F_k)$. For simplicity we assume in this
discussion that
all periodic points of $\ov\alpha $ are fixed points.

Theorem \pt{} asserts
that, as $n\to+\infty$, the sequence
$\da^n(X)$ converges to some $h_\alpha (X)\in\fix\da\inc\bo F_k$ for
every $X\in\bo F_k$.
Let $U$ be the open set $\bo F_k-\fix\da $. Is the function
$h_\alpha $   locally constant on $U$?
 Is the
convergence  of $\da^n$ to $h_\alpha $   locally uniform on $U$? 

Elements of $h_\alpha (U)$ not in $\bo\fix\alpha $ are attracting fixed
points. The action of
$\fix\alpha $ on the set of attracting fixed
points has at most $2k $ orbits [\GJLL].
It is proved in [\Hi] that the action  of $\fix\alpha $ on   $h_\alpha
(U)\cap
\bo\fix\alpha
$ also has finitely many
    orbits.

\subhead \subsect Hyperbolic groups\endsubhead

As shown in \S\kern.15em 10, some of our  results  about automorphisms of free
groups may be extended to   hyperbolic groups.  Another example is
Proposition \bfh{},
which  readily extends to non-elementary, virtually torsion-free,
hyperbolic groups.
    On the other hand, we do not know how to prove that exponentially
growing automorphisms
of an infinitely-ended hyperbolic group $\Gamma  $ have
asymptotically periodic dynamics
on
$\bo\Gamma $.

Let $\alpha \in\Aut(\Gamma )$, with $\Gamma $ hyperbolic. Is there a
bound depending only 
on
$\Gamma $ for the number of ($\fix\alpha
$)-orbits of attracting \fp s of $\da$? Can one associate  a set of growth
rates $\Lambda (\Phi )\inc(1,+\infty)$ to $\Phi \in\Out (\Gamma )$ as
in [\LLdeux]? As
when $\Gamma =F_k$, elements $\lambda \in\Lambda (\Phi )$ should be
the growth rates of
conjugacy classes under iteration of $\Phi $, and also the  rates of
convergence
towards \fp s with respect to the canonical \hr{} structure  on $\bo \Gamma $
(see [\LLdeux]). They should be either dilation factors of pseudo-Anosov
homeomorphisms associated to hanging Fuchsian subgroups, or eigenvalues of
a transition matrix associated to a (relative) train track coming from a
 free product structure. In particular, there should be an upper bound to
the cardinality of
$\Lambda (\Phi )$ that only depends  on $\Gamma $ (not on $\Phi $).

We conclude with the following observation (see [\BHa] for the definition
of $C^*$-simplicity, and [\Be]  for related   results):

\thm{Proposition \sta} Let $\Gamma $ be   a non-elementary torsion-free
hyperbolic group. If
$H\inc\Aut(\Gamma )$ contains all inner automorphisms, then $H$ is
$C^*$-simple.
\fthm

\demo{Proof} We show that the action of $H$ on $\bo \Gamma $ satisfies
the hypotheses of Proposition 1.1 of [\BHa].  Since all non-trivial inner
automorphisms have North-South dynamics, it suffices to show that $\fix\da$
is nowhere dense in $\bo \Gamma $ whenever $\alpha $ is a nontrivial
automorphism of $\Gamma $. 

Assume this is false. Since fixed points of $\da$ not in $\bo\fix\alpha $ are
 isolated (see \S\kern.15em 4 of [\LLdeux]),   $\bo\fix\alpha $ has nonempty
interior. As
$\fix\alpha
$ is quasiconvex [\Ne], it has finite index. It is a general fact (valid
for automorphisms of an arbitrary finitely generated group) that the
centralizer of $\alpha (g)g\mi$ in $\Gamma $ then has finite index for any
$g\in\Gamma
$. Since $ \Gamma $ is torsion-free   and non-elementary, this implies
$\alpha (g)g\mi=1$, so 
$\alpha
$ is the identity. 
\cqfd\enddemo

\bigskip
\Refs
\widestnumber\no{99}
\refno=0

\bref \by E. B\'edos\paper
Discrete groups and simple $C\sp *$-algebras\jour
Math. Proc. Cambridge Philos. Soc. \vol109 \yr1991\pages 521--537
\endref

\bref \by M. Bestvina, M. Feighn \paper   Stable actions of groups
on real trees \jour Invent. Math. \vol121 \yr1995\pages 287--321
\endref

\bref \by M. Bestvina, M. Feighn, M. Handel\paper Laminations, trees, and
irreducible
automorphisms of free groups\jour GAFA\vol7\yr1997\pages215--244
\endref

\bref \by M. Bestvina, M. Feighn, M. Handel\paper The Tits alternative for
$\Out(F_n)$,  I, Dynamics of exponentially-growing
       automorphisms \jour Ann. of Math. \vol 151 \yr2000\pages 517--623
\endref

\bref \by M. Bestvina, M. Feighn, M. Handel
\paper 
The Tits alternative for $\Out(F_n)$. II. A Kolchin type theorem
\jour Ann. of Math. \vol  161  \yr2005\pages 1--59 
\endref

\bref \by M. Bestvina, M. Feighn, M. Handel\paper
Solvable subgroups of $\Out(F_n)$ are virtually abelian \jour
Geometriae Dedicata
\vol 104\pages71--96\yr2004
\endref

\bref \by M. Bestvina, M. Handel\paper Train tracks for automorphisms
of the free group \jour Ann. Math.\vol135 \yr1992\pages1--51 \endref

\bref \by B. Bowditch \paper Cut points and canonical splittings of
hyperbolic groups\jour Acta Math. \vol180\yr1998\pages145--186\endref

\bref\by M.R. Bridson, P. de la Harpe\paper Mapping class groups and outer
automorphism groups of free groups are $C^*$-simple \jour J. Funct. Anal.
\vol212\yr2004\pages 195--205
\endref

\bref \by M.M. Cohen, M. Lustig \paper On the dynamics and the fixed
subgroup of a free group automorphism\jour Inv.
Math.\vol96\yr1989\pages613--638\endref

\bref  \by M.M. Cohen, M. Lustig \paper Very small group actions on
\Rt s and Dehn twist automorphisms\jour Topology\vol34\yr1995\pages575--617
\endref

\bref\by D.J. Collins, E.C. Turner\paper Efficient representatives for
automorphisms of free products\jour Michigan Math. Jour.\vol41\yr1994\pages
443--464\endref

\bref\by F. Dahmani\paper Combination of convergence groups\jour Geom. \&
Top.\vol7\yr2003\pages933--963
\endref

\bref \by W. Dicks, E. Ventura\paper The group fixed by a family of injective
endomorphisms of a free group\jour Contemp. Math.\vol195\yr1996\endref

\bref\by B. Farb \paper Relatively hyperbolic groups\jour
GAFA\vol8\yr1998\pages1--31\endref

\bref\by D. Gaboriau\jour in preparation
\endref

\bref\by D. Gaboriau, A. Jaeger, G. Levitt, M. Lustig\paper An index
for counting
fixed points of automorphisms of free groups\jour Duke Math. Jour.
\vol93\yr1998\pages425--452
\endref

\bref  \by D. Gaboriau, G. Levitt\paper The rank of actions on \Rt s
\jour Ann. Sc. ENS \vol28\yr1995\pages549--570 \endref

\bref\by V. Guirardel\paper Reading small actions
of a one-ended hyperbolic group on $\bold R$-trees from its JSJ splitting
\jour Amer. J. Math. \vol122
       \yr2000\pages 667--688
\endref

\bref\by V. Guirardel\paper 
C\oe ur et nombre d'intersection pour les actions de groupes sur les
arbres\jour  Ann. Sci. \'Ecole Norm. Sup. \vol  38  \yr2005\pages 847--888 
\endref

\bref  \by G.H. Hardy, E.M. Wright\paper An introduction to the
theory of numbers
\jour Oxford Univ. Press \endref

\bref\by A. Hilion\paper Dynamique des automorphismes des groupes libres\jour
Thesis (Toulouse, 2004)\endref

\bref\by I. Kapovich\paper Quasiconvexity and amalgams\jour Int. Jour. Alg.
Comp.\vol7\yr1997\pages771--811\endref

\bref\by I. Kapovich, D. Wise \paper The equivalence of some residual
properties of
word-hyperbolic groups \jour J. Algebra \vol223\pages562--583\yr 2000\endref

\bref\by G. Levitt\paper Automorphisms of hyperbolic groups and graphs of
groups\jour Geometriae Dedicata \vol114 \yr2005\pages 49--70 \endref

\bref\by G. Levitt, M. Lustig\paper Most automorphisms of a
hyperbolic group have 
very simple dynamics \jour Ann. Sc. ENS\vol33\yr2000\pages507--517\endref

\bref\by G. Levitt, M. Lustig\paper Periodic ends, growth rates, H\"older
dynamics for automorphisms of free groups \jour Comm. Math.
Helv.\yr2000\vol75\pages415--430\endref

\bref\by G. Levitt, M. Lustig\paper Irreducible automorphisms of $F_n$
have North-South dynamics on compactified outer space\jour
Jour. Inst. Math. Jussieu\vol2\yr2003\pages 59--72\endref

\bref\by G. Levitt, M. Lustig\paper Growth rates for automorphisms of
free groups \jour
in preparation\endref

\bref\by M. Lustig \paper A discrete action on the product of two
non-simplicial
\Rt s\jour preprint
\endref

\bref \by M. Lustig \paper Structure and conjugacy for automorphisms
of free groups
    I, II
\jour  MPI-Preprint  Series 130 (2000), 4 (2001)
\endref

\bref\by R.T. Miller\paper  Geodesic laminations from Nielsen's viewpoint \jour
Adv. in Math.\yr 1982\vol45\pages189--212
\endref

\bref\by W.D. Neumann\paper The fixed subgroup of an automorphism of a word
hyperbolic group is rational\jour Inv.
Math.\vol110\yr1992\pages147--150\endref

\bref\by F. Paulin  \paper Outer automorphisms of hyperbolic groups and small
actions on
$R$-trees,
{\rm pp\. 331--343 in ``Arboreal group theory (R.C. Alperin, ed.)''}\jour
MSRI Publ. 19\publ Springer Verlag\yr 1991 \endref

\bref \by F. Paulin\paper Sur les automorphismes ext\'erieurs des
groupes hyperboliques \jour Ann. Scient. \'Ec. Norm. Sup.
\vol30\yr1997\pages147--167\endref

\bref \by Z. Sela\paper The Nielsen-Thurston classification and automorphisms
of a free group I\jour Duke Math. Jour.\vol84\yr1996\pages 379--397\endref

\bref\by J. Shor\paper A Scott conjecture for hyperbolic groups\jour
preprint\endref

\bref\by G.A. Swarup
\paper Proof of a weak hyperbolization theorem \jour Q. J. Math. \vol51
    \yr2000\pages529--533
\endref

\endRefs

\address  G.L.:  LMNO, umr cnrs 6139, BP 5186,
Universit\'e de Caen, 14032 Caen Cedex, France.\endaddress\email 
levitt\@math.unicaen.fr{}{}{}{}{}\endemail

\address  M.L.:
Laboratoire de math\'ematiques fondamentales et appliqu\'ees, Universit\'e
d'Aix-Marseille III, 13397 Marseille Cedex 20, France.\endaddress

\email
Martin.Lustig\@math.u-3mrs.fr\endemail

\enddocument